\documentclass[letterpaper,10pt, journal,twoside]{IEEEtran}  % Comment this line out if you need a4paper

\IEEEoverridecommandlockouts                              % This command is only needed if 
\usepackage{amsmath,amssymb,color,cite}%,showkeys}
\usepackage{graphicx}

\usepackage{latexsym}
\usepackage{verbatim}
\usepackage{cite}
\usepackage[english]{babel}
\usepackage[latin1]{inputenc}
\usepackage{enumerate}
\usepackage{amsmath,amsthm}

\usepackage{tikz}
\usepackage{amsfonts}
\usepackage{amssymb}
\usepackage{amsbsy}
\usepackage{bm}

    \let\leq\leqslant
    \let\geq\geqslant

\newcommand{\bmu}{\bm{u}}
\newcommand{\bmx}{\bm{x}}

\newcommand{\bmw}{\bm{w}}
\newcommand{\bmy}{\bm{y}}
\newcommand{\xm}{X_-}
\newcommand{\xp}{X_+}
\newcommand{\um}{U_-}

\newcommand{\xmt}{X_-^\top}
\newcommand{\xpt}{X_+^\top}
\newcommand{\umt}{U_-^\top}

\DeclareMathOperator{\rs}{rs}
\DeclareMathOperator{\iso}{i/s/o}
\DeclareMathOperator{\is}{i/s}
\DeclareMathOperator{\io}{i/o}
\DeclareMathOperator{\im}{im}

\DeclareMathOperator{\rank}{rank}
\DeclareMathOperator{\trace}{tr}

\DeclareMathOperator{\Real}{Re}
\DeclareMathOperator{\Imag}{Im}

\let\leq\leqslant
\let\geq\geqslant

\newcommand{\calD}{\ensuremath{\mathcal{D}}}

\newcommand{\calH}{\ensuremath{\mathcal{H}}}

\newcommand{\calK}{\ensuremath{\mathcal{K}}}
\newcommand{\calL}{\ensuremath{\mathcal{L}}}

\newcommand{\calP}{\ensuremath{\mathcal{P}}}

\newcommand{\calS}{\ensuremath{\mathcal{S}}}

\newcommand{\baru}{\ensuremath{\bar{u}}}

\newcommand{\barP}{\ensuremath{\bar{P}}}

\newcommand{\bmat}{\begin{matrix}}
\newcommand{\emat}{\end{matrix}}
\newcommand{\bbm}{\begin{bmatrix}}
\newcommand{\ebm}{\end{bmatrix}}
\newcommand{\bbma}{\begin{bmatrix*}[r]}
\newcommand{\ebma}{\end{bmatrix*}}
\newcommand{\bpm}{\begin{pmatrix}}
\newcommand{\epm}{\end{pmatrix}}
\newcommand{\bvm}{\begin{vmatrix}}
\newcommand{\evm}{\end{vmatrix}}
\newcommand{\bse}{\begin{subequations}}
\newcommand{\ese}{\end{subequations}}
\newcommand{\beq}{\begin{equation}}
\newcommand{\eeq}{\end{equation}}
\newcommand{\ben}{\renewcommand{\labelenumi}{\arabic{enumi}.}
\renewcommand{\theenumi}{\arabic{enumi}}\begin{enumerate}}

\newcommand{\een}{\end{enumerate}}

\newcommand{\beni}{\renewcommand{\labelenumi}{\roman{enumi}.}
\renewcommand{\theenumi}{\roman{enumi}}\begin{enumerate}}

\newcommand{\eeni}{\end{enumerate}}

\newcommand{\bena}{\renewcommand{\labelenumi}{\alph{enumi}.}
\renewcommand{\theenumi}{\alph{enumi}}\begin{enumerate}}

\newcommand{\eena}{\end{enumerate}}

\newcommand{\bit}{\begin{itemize}}
\newcommand{\eit}{\end{itemize}}
\newcommand{\bthe}{\begin{theorem}}
\newcommand{\ethe}{\end{theorem}}
\newcommand{\blem}{\begin{lemma}}
\newcommand{\elem}{\end{lemma}}
\newcommand{\bprop}{\begin{proposition}}
\newcommand{\eprop}{\end{proposition}}
\newcommand{\bex}{\begin{example}}
\newcommand{\eex}{\end{example}}
\newcommand{\bas}{\begin{assumption}}
\newcommand{\eas}{\end{assumption}}
\newcommand{\bre}{\begin{remark}}
\newcommand{\ere}{\end{remark}}
\newcommand{\bcor}{\begin{corollary}}
\newcommand{\ecor}{\end{corollary}}
\newcommand{\bdfn}{\begin{definition}}
\newcommand{\edfn}{\end{definition}}
\newcommand{\bcon}{\begin{conjecture}}
\newcommand{\econ}{\end{conjecture}}

\newcommand{\half}{\ensuremath{\frac{1}{2}}}
\newcommand{\inv}{\ensuremath{^{-1}}}
\newcommand{\pset}[1]{\ensuremath{\{#1\}}}

\newcommand{\zset}{\ensuremath{\pset{0}}}

\newcommand{\set}[2]{\ensuremath{\{#1\mid #2\}}}

\newcommand{\R}{\ensuremath{\mathbb R}}
\newcommand{\C}{\ensuremath{\mathbb C}}

\newcommand{\BP}{\noindent{\bf Proof. }}
\newcommand{\EP}{\hspace*{\fill} $\blacksquare$\bigskip\noindent}

\newcounter{todocounter}
\usepackage[colorinlistoftodos]{todonotes}

\newtheorem{theorem}{Theorem}
\newtheorem{lemma}[theorem]{Lemma}
\newtheorem{cor}[theorem]{Corollary}
\theoremstyle{definition}
\newtheorem{definition}[theorem]{Definition}
\newtheorem{example}[theorem]{Example}
\newtheorem{prop}[theorem]{Proposition}
\newtheorem{problem}{Problem}
\newtheorem{conjecture}[theorem]{Conjecture}
\newtheorem{remark}[theorem]{Remark}

\title{Data informativity: a new perspective on data-driven analysis and control}

\author{Henk J. van Waarde, Jaap Eising, Harry L. Trentelman, and 
M. Kanat Camlibel
	\thanks{The authors are with the Bernoulli Institute for Mathematics, Computer Science, and Artificial Intelligence, University of Groningen, Nij\-enborgh 9, 9747 AG, Groningen, The Netherlands. Henk van Waarde is also with the Engineering and Technology Institute Groningen, University of Groningen, Nij\-enborgh 4, 9747 AG, Groningen, The Netherlands. (email: {\footnotesize{\tt h.j.van.waarde@rug.nl;j.eising@rug.nl; h.l.trentelman@rug.nl;m.k.camlibel@rug.nl}}).}
}

\begin{document}

\maketitle
\begin{abstract}
The use of persistently exciting data has recently been popularized in the context of data-driven analysis and control. Such data have been used to assess system theoretic properties and to construct control laws, without using a system model. Persistency of excitation is a strong condition that also allows unique identification of the underlying dynamical system from the data within a given model class. In this paper, we develop a new framework in order to work with data that are {\em not\/} necessarily persistently exciting. Within this framework, we investigate necessary and sufficient conditions on the informati\-vity of data for several data-driven analysis and control problems. For certain analysis and design problems, our results reveal that persistency of excitation is not necessary. In fact, in these cases data-driven analysis/control is possible while the combination of (unique) system identification and model-based control is not. For certain other control problems, our results justify the use of persistently exciting data as data-driven control is possible only with data that are informative for system identification.

\end{abstract}

\section{Introduction}

\IEEEPARstart{O}{ne} of the main paradigms in the field of systems and control is that of \emph{model-based} control. Indeed, many control design techniques rely on a system model, represented by e.g. a state-space system or transfer function. In practice, system models are rarely known a priori and have to be identified from measured data using system identification methods such as prediction error \cite{Ljung1999} or subspace identification \cite{vanOverschee1996}. As a consequence, the use of model-based control techniques inherently leads to a two-step control procedure consisting of system identification followed by control design. 

In contrast, \emph{data-driven} control aims to bypass this two-step procedure by constructing controllers directly from data, without (explicitly) identifying a system model. This direct approach is not only attractive from a conceptual point of view but can also be useful in situations where system identification is difficult or even impossible because the data do not give sufficient information.
%
%
%; it also proves to be invaluable in situations where system identification is difficult, for example when the data are not sufficiently informative to uniquely identify a system, or when the data are corrupted by noise. 

The first contribution to data-driven control is often attributed to Ziegler and Nichols for their work on tuning PID controllers \cite{Ziegler1942}. Adaptive control \cite{Astrom1989}, iterative feedback tuning \cite{Hjalmarsson1994,Hjalmarsson1998} and unfalsified control \cite{Safonov1997} can also be regarded as classical data-driven control techniques. More recently, the problem of finding optimal controllers from data has received considerable attention \cite{Bradtke1993,Skelton1994,Furuta1995,Shi1998,Favoreel1999b,Aangenent2005,Markovsky2007,Pang2018,daSilva2019,Baggio2019,Mukherjee2018,Alemzadeh2019}. The proposed solutions to this problem are quite varied, ranging from the use of batch-form Riccati equations \cite{Skelton1994} to approaches that apply reinforcement learning \cite{Bradtke1993}. Additional noteworthy data-driven control pro\-blems include predictive control \cite{Favoreel1999,Salvador2018,Coulson2019}, model reference control \cite{Formentin2013,Campestrini2017} and (intelligent) PID control \cite{Keel2008,Fliess2013}. For more references and classifications of data-driven control techniques, we refer to the survey \cite{Hou2013}. 

In addition to control problems, also \emph{analysis} problems have been studied within a data-based framework. The authors of \cite{Park2009} analyze the stability of an input/output system using time series data. The papers \cite{Wang2011,Liu2014,Niu2017,Zhou2018} deal with data-based controllability and observability analysis. Moreover, the problem of verifying dissipativity on the basis of measured system trajectories has been studied in \cite{Maupong2017b,Romer2017,Berberich2019,Romer2019}.

A result that is becoming increasingly popular in the study of data-driven problems is the so-called \emph{fundamental lemma} by Willems and coworkers \cite{Willems2005}. This result roughly states that all possible trajectories of a linear time-invariant system can be obtained from any given trajectory whose input component is persistently exciting. The fundamental lemma has clear implications for system identification. Indeed, it provides criteria under which the data are sufficiently informative to uniquely identify the system model within a given model class. In addition, the result has also been applied to data-driven control problems. The idea is that control laws can be obtained directly from data, with the underlying mechanism that the system is represented implicitly by the so-called Hankel matrix of a measured trajectory. This framework has led to several interesting control strategies, first in a behavioral setting \cite{Markovsky2007,Markovsky2008,Maupong2017}, and more recently in the context of state-space systems \cite{Coulson2019,DePersis2019,Berberich2019,Romer2019,Berberich2019b,Huang2019}.

The above approaches all use persistently exciting data in the control design, meaning that one could (hypothetically) identify the system model from the same data. An intriguing question is therefore the following: is it possible to obtain a controller from data that are \emph{not} informative enough to uniquely identify the system? An affirmative answer would be remarkable, since it would highlight situations in which data-driven control is more powerful than the combination of system identification and model-based control. On the other hand, a negative answer would also be significant, as it would give a theoretic justification for the use of persistently exciting data for data-driven analysis and control.

To address the above question, this paper introduces a general framework to study data informativity problems for data-driven analysis and control. Specifically, our contributions are the following:

\newpage 

\begin{enumerate}
	\item Inspired by the concept of data informativity in system identification \cite{Ljung1999,Gevers2009,Gevers2013}, we introduce a general notion of informativity for data-driven analysis and control. 
	\item We study the data-driven analysis of several system theoretic properties like stability, stabilizability and controllability. For each of these problems, we provide necessary and sufficient conditions under which the data are informative for this property, i.e., conditions required to ascertain the system's property from data.
	\item We study data-driven control problems such as stabilization by state feedback, stabilization by dynamic measurement feedback, deadbeat control and linear quadratic regulation. In each of the cases, we give conditions under which the data are informative for controller design. 
	\item For each of the studied control problems, we develop methods to compute a controller from data, assuming that the informativity conditions are satisfied. 
\end{enumerate}

Our work has multiple noteworthy implications. First of all, we show that for problems like stabilization by state feedback, the corresponding informativity conditions on the data are \emph{weaker} than those for system identification. This implies that a stabilizing feedback can be obtained from data that are not sufficiently informative to uniquely identify the system. 

Moreover, for problems such as linear quadratic regu\-lation (LQR), we show that the informativity conditions are essentially the same as for system identification. Therefore, our results provide a theoretic justification for imposing the strong persistency of excitation conditions in prior work on the LQR problem, such as \cite{Markovsky2007} and \cite{DePersis2019}.

The paper is organized as follows. In Section \ref{sec:prob} we introduce the problem at a conceptual level. Subsequently, in Section \ref{sec:dd analysis} we provide data informativity conditions for controllability and stabilizability. Section \ref{sec:no output} deals with data-driven control problems with input/state data. Next, Section \ref{sec:output} discusses control problems where ouput data plays a role. Finally, Section~\ref{sec:conc} contains our conclusions and suggestions for future work. 

\section{Problem formulation}\label{sec:prob}
In this section we will first introduce the \textit{informativity framework} for data-driven analysis and control in a fairly abstract manner.

Let $\Sigma$ be a model class, i.e. a given set of systems containing the `true' system denoted by $\calS$. We assume that the `true' system $\calS$ is not known but that we have access to a set of data, $\calD$, which are generated by this system. In this paper we are interested in assessing system-theoretic properties of $\calS$ and designing control laws for it from the data $\calD$.

Given the data $\calD$, we define $\Sigma_\calD\subseteq\Sigma$ to be the set of all systems that are consistent with the data $\calD$, i.e. that could also have generated these data.

We first focus on data-driven analysis. Let $\calP$ be a system-theoretic property. We will denote the set of all systems within $\Sigma$ having this property by $\Sigma_\calP$.

Now suppose we are interested in the question whether our `true' system $\calS$ has the property $\calP$. As the only information we have to base our answer on are the data $\calD$ obtained from the system, we can only conclude that the `true' system has property $\calP$ if \textit{all} systems consistent with the data $\calD$ have the property $\calP$. This leads to the following definition:

\begin{definition}[Informativity]\label{def:informativity}
	We say that the data $\calD$ are \textit{informative} for property $\calP$ if $\Sigma_\calD \subseteq\Sigma_\calP$. 
\end{definition}

Next, we illustrate the above abstract setup by an example.

\begin{example}\label{e:exmp1} 
For given $n$ and $m$, let $\Sigma$ be the set of all discrete-time linear input/state systems of the form
$$
\bmx(t+1) = A\bmx(t) + B\bmu(t)
$$
where $\bmx$ is the $n$-dimensional state and $\bmu$ is the $m$-dimensional input. Let the `true' system $\calS$ be represented by the matrices $(A_s,B_s)$. 

An example of a data set $\calD$ arises when considering data-driven problems on the basis of input and state measurements. Suppose that we collect input/state data on $q$ time intervals $\pset{0,1,\ldots,T_i}$ for $i=1,2,\ldots,q$. Let
\begin{subequations}\label{eq: UXdata}
\begin{align}
U^i_-& := \bbm u^i(0) & u^i(1) & \cdots & u^i(T_i-1)\ebm, \\
X^i& := \bbm x^i(0) & x^i(1) & \cdots & x^i(T_i)\ebm \label{eq: UXdata2}
\end{align}
\end{subequations}
denote the input and state data on the $i$-th interval. By defining
\begin{subequations}\label{eq: def of X- X+}
\begin{align}
X^i_-& := \bbm x^i(0) & x^i(1) & \cdots & x^i(T_i-1) \ebm, \\
X^i_+& := \bbm x^i(1) & x^i(2) & \cdots & x^i(T_i) \ebm,
\end{align}	
\end{subequations}
we clearly have $X^i_+=A_sX^i_-+B_sU^i_-$ for each $i$ because the `true' system is assumed to generate the data. Now, introduce the notation
\bse\label{eq: UXdatanew}
\begin{alignat}{3}
U_-&:=\bbm U^1_-&\cdots& U^q_-\ebm,&\quad X&:=\bbm X^1&\cdots& X^q\ebm,\label{eq: UXdatanew1}\\
X_-&:=\bbm X^1_-&\cdots& X^q_-\ebm,&\quad X_+&:=\bbm X^1_+&\cdots& X^q_+\ebm.\label{eq: UXdatanew2}
\end{alignat}
\ese
We then define the data as $\calD:=(U_-,X)$. In this case, the set $\Sigma_\calD$ is equal to $\Sigma_{(U_-,X)}$ defined by 
	\begin{equation}\label{eq: SigmaD}
	\Sigma_{(U_-,X)} := \left\{ (A,B) \mid X_+= \bbm A&B \ebm
	\begin{bmatrix}
	X_-\\U_-
	\end{bmatrix} \right\}.
	\end{equation}
Clearly, we have $(A_s,B_s)\in\Sigma_\calD$.

Suppose that we are interested in the system-theoretic property $\calP$ of \textit{stabilizability}. The corresponding set $\Sigma_\calP$ is then equal to $\Sigma_{\mathrm{stab}}$ defined by
\[ \Sigma_{\mathrm{stab}} :=\{ (A,B) \mid (A,B) \text{ is stabilizable}\}. \] 
Then, the data $(U_-,X)$ are informative for stabilizability if $\Sigma_{(U_-,X)} \subseteq\Sigma_{\mathrm{stab}}$. That is, if all systems consistent with the input/state measurements are stabilizable. 
\end{example}

In general, if the `true' system $\calS$ can be uniquely determined from the data $\calD$, that is $\Sigma_\calD=\pset{\calS}$ and $\calS$ has the property $\calP$, then it is evident that the data $\calD$ are informative for $\calP$. However, the converse may not be true: $\Sigma_\calD$ might contain many systems, all of which have property $\calP$. This paper is interested in necessary \textit{and} sufficient conditions for informativity of the data. Such conditions reveal the minimal amount of information required to assess the property $\calP$. A natural problem statement is therefore the following:

\begin{problem}[Informativity problem]\label{prob:general}
	Provide necessary and sufficient conditions on $\calD$ under which the data are informative for property $\calP$. 
\end{problem}

The above gives us a general framework to deal with data-driven analysis problems. Such analysis problems will be the main focus of Section~\ref{sec:dd analysis}. 

This paper also deals with data-driven control problems. The objective in such problems is the data-based design of controllers such that the closed loop system, obtained from the interconnection of the `true' system $\calS$ and the controller, has a specified property. 

As for the analysis problem, we have only the information from the data to base our design on. Therefore, we can only guarantee our control objective if the designed controller imposes the specified property when interconnected with \textit{any} system from the set $\Sigma_\calD$. 

For the framework to allow for data-driven control problems, we will consider a system-theoretic property $\calP(\calK)$ that depends on a given controller $\calK$. For properties such as these, we have the following variant of informativity: 

\begin{definition}[Informativity for control]\label{def:par informativity}
	We say that the data $\calD$ are \textit{informative} for the property $\calP(\cdot)$ if there exists a controller $\calK$ such that $\Sigma_\calD \subseteq\Sigma_{\calP(\calK)}$. 
\end{definition}

\begin{example} 
	For systems and data like in Example~\ref{e:exmp1}, we can take the controller $\calK=K\in\mathbb{R}^{m\times n}$ and the property $\calP(\calK):$ `interconnection with the state feedback $K$ yields a stable closed loop system'. The corresponding set of systems $\Sigma_{\calP(\calK)}$ is equal to $\Sigma_K$ defined by
	\[ \Sigma_K = \{(A,B) \mid A+BK \text{ is stable}\footnote{We say that a matrix is {\em stable\/} if all its eigenvalues are contained in the open unit disk.} \}. \] 
\end{example}

The first step in any data-driven control problem is to determine whether it is possible to obtain a suitable controller from given data. This leads to the following informativity problem:

\begin{problem}[Informativity problem for control]\label{prob:parametrized}
Provide necessary and sufficient conditions on $\calD$ under which there exists a controller $\calK$ such that the data are informative for property $\calP(\calK)$. 
\end{problem}

The second step of data-driven control involves the design of a suitable controller. In terms of our framework, this can be stated as:

\begin{problem}[Control design problem]\label{prob:design}
	Under the assumption that the data $\calD$ are informative for property $\calP(\cdot)$, find a controller $\calK$ such that $\Sigma_\calD \subseteq \Sigma_{\calP(\calK)}$. 
\end{problem}

As stated in the introduction, we will highlight the strength of this framework by solving multiple problems. We stress that throughout the paper it is assumed that the data are \textit{given} and are \textit{not corrupted by noise}.

\section{Data-driven analysis}\label{sec:dd analysis} 
In this section, we will study data-driven analysis of controllability and stabilizability given input and state measurements. As in Example~\ref{e:exmp1}, consider the discrete-time linear system
\begin{equation}\label{e:disc}
\bmx(t+1)=A_s\bmx(t)+B_s\bmu(t).
\end{equation}
We will consider data consisting of input and state measurements. We define the matrices $U_-$ and $X$ as in \eqref{eq: UXdatanew1} and define $X_-$ and $X_+$ as in \eqref{eq: UXdatanew2}. The set of all systems compatible with these data was introduced in \eqref{eq: SigmaD}. In order to stress that we deal with input/state data, we rename it here as
\begin{equation}\label{eq:sigma is} \Sigma_{\is} := \left\{ (A,B) \mid X_+= \bbm A&B \ebm 
\begin{bmatrix}
X_-\\U_-
\end{bmatrix} \right\}. \end{equation}
Note that the defining equation of \eqref{eq:sigma is} is a system of linear equations in the unknowns $A$ and $B$. The solution space of the corresponding homogeneous equations is denoted by $\Sigma_{\is}^0$ and is equal to
\begin{equation}\label{eq:sigma is 0}
\Sigma_{\is}^0 := \left\{(A_0,B_0) \mid 0 = \bbm A_0 & B_0 \ebm
\begin{bmatrix}
X_- \\ U_-
\end{bmatrix}\right\}.
\end{equation}

We consider the problem of data-driven analysis for systems of the form \eqref{e:disc}. If $(A_s,B_s)$ is the only system that explains the data, data-driven analysis could be performed by first identifying this system and then analyzing its properties. It is therefore of interest to know under which conditions there is only one system that explains the data. 

\begin{definition}
	We say that the data $(U_-,X)$ are {\em informative for system identification} if $\Sigma_{\is} = \{(A_s,B_s)\}$. 
\end{definition}

It is straightforward to derive the following result:

\begin{prop}\label{prop:sys ident}
	The data $(U_-,X)$ are informative for system identification if and only if 
	\begin{equation}\label{eq:inf for sys ident} \rank \begin{bmatrix} X_- \\ U_- \end{bmatrix} = n+m. \end{equation}
	Furthermore, if \eqref{eq:inf for sys ident} holds, there exists a right inverse\footnote{Note that $\begin{bmatrix}
V_1 & V_2
\end{bmatrix}$ is not unique whenever $T > n+m$.} $\begin{bmatrix}	V_1 & V_2 \end{bmatrix}$ such that 
	\begin{equation}\label{eq:V1 V2} \begin{bmatrix} X_- \\ U_- \end{bmatrix} \begin{bmatrix} 	V_1 & V_2	\end{bmatrix} = \begin{bmatrix} I & 0 \\ 0& I\end{bmatrix}, \end{equation}
	and for any such right inverse $A_s=X_+ V_1$ and $B_s= X_+ V_2$. 
\end{prop}
As we will show in this section, the condition \eqref{eq:inf for sys ident} is not necessary for data-driven analysis in general. We now proceed by studying data-driven analysis of controllability and stabilizability. Recall the Hautus test \cite[Theorem 3.13]{Trentelman2001} for controllability: a system $(A,B)$ is controllable if and only if
\begin{equation}\label{eq:Hautus} \rank \begin{bmatrix} A-\lambda I & B \end{bmatrix} =n \end{equation} 
for all $\lambda\in \mathbb{C}$. For stabilizability, we require that \eqref{eq:Hautus} holds for all $\lambda$ outside the open unit disc.

Now, we introduce the following sets of systems:
\begin{align*}
\Sigma_{\text{cont}} &:= \{ (A,B) \mid (A,B) \text{ is controllable} \} \\
\Sigma_{\text{stab}} &:= \{ (A,B) \mid (A,B) \text{ is stabilizable} \}.
\end{align*} 

Using Definition~\ref{def:informativity}, we obtain the notions of \textit{informativity for controllability} and \textit{stabilizability}. To be precise:

\begin{definition}
	We say that the data $(U_-,X)$ are {\em informative for controllability} if $\Sigma_{\is} \subseteq \Sigma_{\text{cont}}$ and {\em informative for stabilizability} if $\Sigma_{\is} \subseteq \Sigma_{\text{stab}}$. 
\end{definition}

In the following theorem, we give necessary and sufficient conditions for the above notions of informativity. The result is remarkable as only data matrices are used to assess controllability and stabilizability. 

\begin{theorem}[Data-driven Hautus tests]
	\label{t:contstab}
	The data $(U_-,X)$ are informative for controllability if and only if 
	\begin{equation}
	\label{eq:rank for cont}
	\rank ( X_+ -\lambda X_-) = n \quad \forall \lambda \in \mathbb{C}.
	\end{equation}
	Similarly, the data $(U_-,X)$ are informative for stabilizability if and only if 
	\begin{equation}\label{eq:rank for stab}
	\rank ( X_+ -\lambda X_-) = n \quad \forall \lambda \in \mathbb{C}\textrm{ with } |\lambda|\geq 1.
	\end{equation}
\end{theorem}
Before proving the theorem, we will discuss some of its implications. We begin with computational issues. 

\begin{remark} Similar to the classical Hautus test, \eqref{eq:rank for cont} and \eqref{eq:rank for stab} can be verified by checking the rank for finitely many complex numbers $\lambda$. Indeed, \eqref{eq:rank for cont} is equivalent to $\rank(X_+)=n$ and 
$$
\rank ( X_+ -\lambda X_-) = n
$$
for all $\lambda\neq 0$ with $\lambda\inv\in\sigma(X_-X_+^\dagger)$, where $X_+^\dagger$ is any right inverse of $X_+$. Here, $\sigma(M)$ denotes the spectrum, i.e. set of eigenvalues of the matrix $M$. Similarly, \eqref{eq:rank for stab} is equivalent to $\rank(X_+-X_-)=n$ and 
$$
\rank ( X_+ -\lambda X_-) = n 
$$
for all $\lambda\neq 1$ with $(\lambda-1)\inv\in\sigma(X_-(X_+-X_-)^\dagger)$, where $(X_+-X_-)^\dagger$ is any right inverse of $X_+-X_-$.

%	In the special case that the data are informative for system identification, we can use the matrix $\begin{bmatrix}
%	V_1 & V_2\end{bmatrix}$ from Proposition~\ref{prop:sys ident} to note that
%	\begin{align*} (X_+-\lambda X_-)\begin{bmatrix}	V_1 & V_2 \end{bmatrix} &= \begin{bmatrix} X_+V_1-\lambda I & X_+V_2\end{bmatrix}\\
%	&= \begin{bmatrix} A_s-\lambda I & B_s\end{bmatrix}. \end{align*}
%	Because $\begin{bmatrix}V_1 & V_2\end{bmatrix}$ has full column rank, the data-driven Hautus test \eqref{eq:rank for cont} then boils down to the Hautus test. In general, however, condition \eqref{eq:rank for cont} is a weaker condition than \eqref{eq:inf for sys ident}. This means that there are situations in which we can conclude controllability from data without being able to identify the `true' system uniquely. 
\end{remark}

A noteworthy point to mention is that there are situations in which we can conclude controllability/stabilizability from the data without being able to identify the `true' system uniquely, as illustrated next.
	
\begin{example}
Suppose that $n=2$, $m=1$, $q=1$, $T_1=2$ and we obtain the data
	\[ X= \begin{bmatrix} 0 & 1 & 0 \\ 0 & 0 & 1 \end{bmatrix} \textrm{ and } U_-= \begin{bmatrix} 1 & 0 \end{bmatrix}. \] 
	This implies that 
	\[X_+ = \begin{bmatrix} 1 & 0 \\ 0 & 1 \end{bmatrix} \textrm{ and } X_- = \begin{bmatrix} 0 & 1 \\ 0& 0\end{bmatrix}. \]
	Clearly, by Theorem~\ref{t:contstab} we see that these data are informative for controllability, as 
	\[ \rank \begin{bmatrix} 1 & -\lambda \\ 0 & 1 \end{bmatrix}= 2 \quad \forall \lambda \in \mathbb{C}.\] 
	As therefore all systems explaining the data are controllable, we conclude that the `true' system is controllable. It is worthwhile to note that the data are not informative for system identification, as 
	\begin{equation}\label{eq:structural} \Sigma_{\is} = \left\lbrace \left( \begin{bmatrix} 0 & a_1 \\ 1 & a_2\end{bmatrix} ,\begin{bmatrix} 1 \\ 0 \end{bmatrix}\right) \mid a_1,a_2 \in\mathbb{R} \right\rbrace . \end{equation}
\end{example}
%\begin{remark}
%    {\color{blue}The data-driven controllability problem is comparable to strong structural controllability (see, e.g., \cite{Jia2019}) in the sense that both problems are concerned with the study of controllability \emph{for all} systems within a given set. What makes these problems different is the particular set of systems. In strong structural controllability, this set originates from a \emph{structure} on the system matrices \cite{Jia2019}, while in data-driven controllability we study an affine set of systems originating from the restrictions imposed by the data.}
%\end{remark} 

\noindent{\bf Proof of Theorem~\ref{t:contstab}.} We will only prove the characterization of informativity for controllability. The proof for stabilizability uses very similar arguments, and is hence omitted. 

Note that the condition \eqref{eq:rank for cont} is equivalent to the implication:
\begin{equation}\label{impcont}
	z \in \mathbb{C}^n, \lambda \in \mathbb{C} \text{ and } z^* X_+ = \lambda z^* X_- \implies z = 0.
\end{equation} 
Suppose that the implication \eqref{impcont} holds. Let $(A,B) \in \Sigma_{\is}$ and suppose that $z^*\begin{bmatrix}
A-\lambda I & B
\end{bmatrix} = 0$. We want to prove that $z = 0$. Note that $z^*\begin{bmatrix}
A-\lambda I & B
\end{bmatrix} = 0$ implies that 
$$
z^* \begin{bmatrix}
A - \lambda I & B
\end{bmatrix} \begin{bmatrix}
X_- \\ U_-
\end{bmatrix} = 0,
$$
or equivalently $z^* X_+ = \lambda z^* X_-$. This means that $z = 0$ by \eqref{impcont}. We conclude that $(A,B)$ is controllable, i.e., $(U_-,X)$ are informative for controllability.

Conversely, suppose that $(U_-,X)$ are informative for controllability. Let $z \in \mathbb{C}^n$ and $\lambda \in \mathbb{C}$ be such that $z^* X_+ = \lambda z^* X_-$. This implies that for all $(A,B) \in \Sigma_{\is}$, we have $z^* \begin{bmatrix}
A & B
\end{bmatrix} \begin{bmatrix}
X_- \\ U_-
\end{bmatrix} = \lambda z^* X_-$. In other words, 

\begin{equation}
\label{leftkernel}
z^* \begin{bmatrix}
A - \lambda I & B
\end{bmatrix} \begin{bmatrix}
X_- \\ U_-
\end{bmatrix} = 0.
\end{equation}

We now distinguish two cases, namely the case that $\lambda$ is real, and the case that $\lambda$ is complex. First suppose that $\lambda$ is real. Without loss of generality, $z$ is real. We want to prove that $z = 0$. Suppose on the contrary that $z \neq 0$ and $z^\top z = 1$. We define the (real) matrices 
$$
\bar{A} := A - z z^\top (A - \lambda I) \text{ and } \bar{B} := B - z z^\top B.
$$
In view of \eqref{leftkernel}, we find that $(\bar{A}, \bar{B}) \in \Sigma_{\is}$. Moreover, 
$$
z^\top \bar{A} = z^\top A - z^\top (A-\lambda I) = \lambda z^\top 
$$
and 
$$
z^\top \bar{B} = z^\top B - z^\top B = 0.
$$
This means that 
$$
z^\top \begin{bmatrix}
\bar{A} - \lambda I & \bar{B}
\end{bmatrix} = 0.
$$
However, this is a contradiction as $(\bar{A},\bar{B})$ is controllable by the hypothesis that $(U_-,X)$ are informative for controllability. We conclude that $z = 0$ which shows that \eqref{impcont} holds for the case that $\lambda$ is real.

Secondly, consider the case that $\lambda$ is complex. We write $z$ as $z = p + iq$, where $p,q \in \mathbb{R}^n$ and $i$ denotes the imaginary unit. If $p$ and $q$ are linearly dependent, then $p = \alpha q$ or $q = \beta p$ for $\alpha,\beta \in \mathbb{R}$. If $p = \alpha q$ then substitution of $z = (\alpha + i)q$ into $z^* X_+ = \lambda z^* X_-$ yields 
\begin{align*}
(\alpha - i)q^\top X_+ &= \lambda (\alpha - i)q^\top X_-,
\end{align*} 
that is, $q^\top X_+ = \lambda q^\top X_-$. As $q^\top X_+$ is real and $\lambda$ is complex, we must have $q^\top X_+ = 0$ and $q^\top X_- = 0$. This means that $z^* X_+ = z^* X_- = 0$, hence $z^* X_+ = \mu z^* X_-$ for any real $\mu$, which means that $z = 0$ by case 1. Using the same arguments, we can show that $z = 0$ if $q = \beta p$. 

It suffices to prove now that $p$ and $q$ are linearly dependent. Suppose on the contrary that $p$ and $q$ are linearly independent. Since $\lambda$ is complex, $n \geq 2$. Therefore, by linear independence of $p$ and $q$ there exist $\eta, \zeta \in \mathbb{R}^n$ such that 
$$
\begin{bmatrix}
p^\top\: \\ q^\top\:
\end{bmatrix}
\begin{bmatrix}
\eta & \zeta 
\end{bmatrix}
= \begin{bmatrix}
1 & 0 \\ 0 & -1
\end{bmatrix}.
$$
We now define the real matrices $\bar{A}$ and $\bar{B}$ as
$$
\begin{bmatrix}\bar{A} & \bar{B}\end{bmatrix} := \begin{bmatrix}
A & B
\end{bmatrix} - \begin{bmatrix}
\eta & \zeta 
\end{bmatrix} \begin{bmatrix}
\Real\left( z^*\begin{bmatrix}
A - \lambda I & B
\end{bmatrix} \right) \\
\Imag\left( z^*\begin{bmatrix}
A - \lambda I & B
\end{bmatrix} \right)
\end{bmatrix}.
$$
By \eqref{leftkernel} we have $(\bar{A},\bar{B}) \in \Sigma_{\is}$. Next, we compute 
\begin{align*}
z^* \begin{bmatrix}\bar{A} & \bar{B}\end{bmatrix} &= z^* \begin{bmatrix}
A & B
\end{bmatrix} - \begin{bmatrix}
1 & i 
\end{bmatrix} \begin{bmatrix}
\Real\left( z^*\begin{bmatrix}
A - \lambda I & B
\end{bmatrix} \right) \\[1mm]
\Imag\left( z^*\begin{bmatrix}
A - \lambda I & B
\end{bmatrix} \right)
\end{bmatrix} \\
&= z^* \begin{bmatrix}
A & B
\end{bmatrix} - z^*\begin{bmatrix}
A - \lambda I & B
\end{bmatrix} \\
&= z^* \begin{bmatrix}
\lambda I & 0
\end{bmatrix}.
\end{align*} 
This implies that $z^* \begin{bmatrix}\bar{A} - \lambda I & \bar{B}\end{bmatrix} = 0$. Using the fact that $(\bar{A},\bar{B})$ is controllable, we conclude that $z = 0$. This is a contradiction with the fact that $p$ and $q$ are linearly independent. Thus $p$ and $q$ are linearly dependent and therefore implication \eqref{impcont} holds. This proves the theorem. 
\EP 

In addition to controllability and stabilizability, we can also study the \emph{stability} of an autonomous system of the form
\begin{equation}\label{e:discaut}
\bmx(t+1)=A_s\bmx(t).
\end{equation}
To this end, let $X$ denote the matrix of state measurements obtained from \eqref{e:discaut}, as defined in \eqref{eq: UXdatanew1}. The set of all autonomous systems compatible with these data is
\[ \Sigma_{\mathrm{s}} :=\left\{ A \mid X_+= AX_- \right\}.\]

Then, we say the data $X$ are \emph{informative for stability} if any matrix $A\in\Sigma_{\text{s}}$ is  stable, i.e. Schur. Using Theorem \ref{t:contstab} we can show that stability can only be concluded if the `true' system can be uniquely identified.

\begin{cor} 
The data $X$ are informative for stability if and only if $X_-$ has full row rank and $X_+X_-^\dagger$ is stable for any right inverse $X_-^\dagger$, equivalently $\Sigma_{\mathrm{s}}=\pset{A_{\mathrm{s}}}$ and $A_{\mathrm{s}}=X_+X_-^\dagger$ is stable.

\end{cor} 

\BP Since the `if' part is evident, we only prove the `only if' part. By taking $B = 0$, it follows from Theorem~\ref{t:contstab} that the data $X$ are informative for stability if and only if 
\begin{equation}\label{eq:rank for stability} \rank ( X_+ -\lambda X_-) = n \quad \forall \lambda \in \mathbb{C}\textrm{ with } |\lambda|\geq 1. \end{equation}
Let $z$ be such that $z^\top X_-=0$. Take $A\in\Sigma_{\mathrm{s}}$ and $\lambda>1$ such that $\lambda$ is not an eigenvalue of $A$. Note that
$$
z^\top(A-\lambda I)\inv  ( X_+ -\lambda X_-)=z^\top X_-=0.
$$
Since $\rank ( X_+ -\lambda X_-) = n$, we may conclude that $z=0$. Hence, $X_-$ has full row rank. Therefore, $\Sigma_{\mathrm{s}}=\pset{A_{\mathrm{s}}}$ where $A_{\mathrm{s}}=X_+X_-^\dagger$ for any right inverse $X_-^\dagger$ and $A_{\mathrm{s}}$ is stable.
\EP

Note that there is a subtle but important difference between the characterizations \eqref{eq:rank for stab} and \eqref{eq:rank for stability}. For the first the data $X$ are assumed to be generated by a system with inputs, whereas the data for the second characterization are generated by an autonomous system. 

\section{Control using input and state data}\label{sec:no output}

In this section we will consider various state feedback control problems on the basis of input/state measurements. First, we will consider the problem of data-driven stabilization by static state feedback, where the data consist of input and state measurements. As described in the problem statement we will look at the informativity and design problems separately as special cases of Problem~\ref{prob:parametrized} and Problem~\ref{prob:design}. We will then use similar techniques to obtain a result for \textit{deadbeat control}. 

After this, we will shift towards the linear quadratic regulator problem, where we wish to find a stabilizing feedback that additionally minimizes a specified quadratic cost. 

\subsection{Stabilization by state feedback} 

In what follows, we will consider the problem of finding a stabilizing controller for the system \eqref{e:disc}, using only the data $(U_-,X)$. To this end, we define the set of systems $(A,B)$ that are stabilized by a given $K$:
$$
\Sigma_K := \{(A,B) \mid A+BK \text{ is stable} \}.
$$

In addition, recall the set $\Sigma_{\is}$ as defined in \eqref{eq:sigma is} and $\Sigma_{\is}^0$ from \eqref{eq:sigma is 0}. In line with Definition~\ref{def:par informativity} we obtain the following notion of informativity for stabilization by state feedback.

\begin{definition}
	We say that the data $(U_-,X)$ are {\em informative for stabilization by state feedback} if there exists a feedback gain $K$ such that $\Sigma_{\is} \subseteq \Sigma_K$.
\end{definition}

\begin{remark} 
At this point, one may wonder about the relation between informativity for stabilizability (as in Section~\ref{sec:dd analysis}) and informativity for stabilization. It is clear that $(U_-,X)$ are informative for stabilizability if $(U_-,X)$ are informative for stabilization by state feedback. However, the reverse statement does not hold in general. This is due to the fact that all systems $(A,B)$ in $\Sigma_{\is}$ may be stabilizable, but there may not be a \emph{common} feedback gain $K$ such that $A+BK$ is stable for all of these systems. Note that the existence of a common stabilizing $K$ for all systems in $\Sigma_{\is}$ is essential, since there is no way to distinguish between the systems in $\Sigma_{\is}$ based on the given data $(U_-,X)$. 
\end{remark}

The following example further illustrates the difference between informativity for stabilizability and informativity for stabilization.

\begin{example}
	Consider the scalar system 
	$$
	\bmx(t+1)= \bmu(t),
	$$
	where $\bmx, \bmu \in \mathbb{R}$. Suppose that $q=1$, $T_1 = 1$ and $x(0) = 0$, $u(0) = 1$ and $x(1) = 1$. This means that $U_- = \begin{bmatrix}1 \end{bmatrix}$ and $X = \begin{bmatrix}0 & 1	\end{bmatrix}$. It can be shown that $\Sigma_{\is} = \{(a,1) \mid a \in \mathbb{R}\}$. Clearly, all systems in $\Sigma_{\is}$ are stabilizable, i.e., $\Sigma_{\is} \subseteq \Sigma_{\text{stab}}$. Nonetheless, the data are not informative for \emph{stabilization}. This is because the systems $(-1,1)$ and $(1,1)$ in $\Sigma_{\is}$ cannot be stabilized by the \emph{same} controller of the form $u(t) = K x(t)$. 
We conclude that informativity of the data for stabilizability does not imply informativity for stabilization by state feedback.
\end{example}

The notion of informativity for stabilization by state feedback is a specific example of informativity for control. As described in Problem~\ref{prob:parametrized}, we will first find necessary and sufficient conditions for informativity for stabilization by state feedback. After this, we will design a corresponding controller, as described in Problem~\ref{prob:design}. 

In order to be able to characterize informativity for stabilization, we first state the following lemma. 

\begin{lemma}
	\label{lemmaF0=0}
	Suppose that the data $(U_-,X)$ are informative for stabilization by state feedback, and let $K$ be a  feedback gain such that $\Sigma_{\is} \subseteq \Sigma_K$. Then $A_0 + B_0 K = 0$ for all $(A_0,B_0) \in \Sigma_{\is}^0$. Equivalently, 
$$
	\im \begin{bmatrix}
I \\ K
\end{bmatrix} \subseteq 
\im \begin{bmatrix}
X_- \\ U_-
\end{bmatrix}.
$$
\end{lemma}

\BP
We first prove that $A_0 + B_0 K$ is \emph{nilpotent} for all $(A_0,B_0) \in \Sigma_{\is}^0$. By hypothesis, $A + BK$ is stable for all $(A,B)\in\Sigma_{\is}$. Let $(A,B) \in \Sigma_{\is}$ and $(A_0,B_0) \in \Sigma_{\is}^0$ and define the matrices $F := A + BK$ and $F_0 := A_0 + B_0 K$. Then, the matrix $F + \alpha F_0$ is stable for all $\alpha \geq 0$. By dividing by $\alpha$, it follows that, for all $\alpha\geq 1$, the spectral radius of the matrix
$$
M_\alpha := \frac{1}{\alpha} F + F_0
$$
is smaller than $1/\alpha$. From the continuity of the spectral radius by taking the limit as $\alpha$ tends to infinity, we see that $F_0=A_0 + B_0 K$ is nilpotent for all $(A_0,B_0) \in \Sigma_{\is}^0$. Note that we have 
$$((A_0 + B_0 K)^TA_0,(A_0 + B_0 K)^TB_0)\in\Sigma_{\is}^0$$
whenever $(A_0,B_0) \in \Sigma_{\is}^0$. This means that $(A_0 + B_0 K)^T(A_0 + B_0 K)$ is nilpotent. Since the only symmetric nilpotent matrix is the zero matrix, we see that $A_0 + B_0 K=0$ for all $(A_0,B_0) \in \Sigma_{\is}^0$. This is equivalent to
$$
\ker \begin{bmatrix}
X_-^\top & U_-^\top 
\end{bmatrix} \subseteq \ker \begin{bmatrix}
I & K^\top 
\end{bmatrix}
$$
which is equivalent to
$
\im \begin{bmatrix}
I \\ K
\end{bmatrix} \subseteq 
\im \begin{bmatrix}
X_- \\ U_-
\end{bmatrix}.
$\EP

%It remains to be shown that $A_0 + B_0 K = 0$ for all $(A_0,B_0) \in \Sigma_{\is}^0$. To prove this, let $S \in \mathbb{R}^{(n+m) \times r}$ be a full column rank matrix such that 
%$$
%\im S = \ker \begin{bmatrix}
%X_- \\ U_- 
%\end{bmatrix}^\top.
%$$
%Clearly, $(A_0,B_0) \in \Sigma_{\is}^0$ if and only if $\begin{bmatrix}
%A_0 & B_0
%\end{bmatrix} = MS^\top$ for some $M \in \mathbb{R}^{n \times r}$. Now, since $A_0 + B_0 K$ is nilpotent for all $(A_0,B_0) \in \Sigma_{\is}^0$, the matrix 
%$$MS^\top \begin{bmatrix} I \\ K \end{bmatrix}$$ 
%is nilpotent for all $M \in \mathbb{R}^{n \times r}$. In particular, the choice $M = \begin{bmatrix}
%I & K^\top 
%\end{bmatrix}S$ shows that  
%\begin{equation}
%\label{symmatrix}
%\begin{bmatrix} I & K \end{bmatrix}^\top S S^\top \begin{bmatrix} I \\ K \end{bmatrix}
%\end{equation}
%is nilpotent. However, the only symmetric nilpotent matrix is the zero matrix. Hence, the matrix in \eqref{symmatrix} is zero and we conclude that 
%$$S^\top \begin{bmatrix}
%I \\ K
%\end{bmatrix} = 0.
%$$
%Therefore, for each $(A_0,B_0) \in \Sigma_{\is}^0$ we have
%$$
%A_0 + B_0 K = MS^\top \begin{bmatrix}
%I \\ K
%\end{bmatrix} = 0,
%$$
%which completes the proof.
%\EP

The previous lemma is instrumental in proving the following theorem that gives necessary and sufficient conditions for informativity for stabilization by state feedback. 

\begin{theorem}
	\label{t:algstab}
	The data $(U_-,X)$ are informative for stabilization by state feedback if and only if the matrix $X_-$ has full row rank and there exists a right inverse $X_-^\dagger$ of $X_-$ such that $X_+ X_-^\dagger$ is stable. 
	
	Moreover, $K$ is such that $\Sigma_{\is} \subseteq \Sigma_K$ if and only if $K = U_- X_-^\dagger$, where $X_-^\dagger$ satisfies the above properties.
\end{theorem}

\BP
To prove the `if' part of the first statement, suppose that $X_-$ has full row rank and there exists a right inverse $X_-^\dagger$ of $X_-$ such that $X_+ X_-^\dagger$ is stable. We define $K := U_- X_-^\dagger$. Next, we see that
\begin{equation}
\label{closedloop}
X_+ X_-^\dagger = \begin{bmatrix}
A & B
\end{bmatrix} \begin{bmatrix}
X_- \\ U_-
\end{bmatrix} X_-^\dagger = A + BK,
\end{equation}
for all $(A,B) \in \Sigma_{\is}$. Therefore, $A + BK$ is stable for all $(A,B) \in \Sigma_{\is}$, i.e., $\Sigma_{\is} \subseteq \Sigma_K$. We conclude that the data $(U_-,X)$ are informative for stabilization by state feedback, proving the `if' part of the first statement. Since $K = U_- X_-^\dagger$ is such that $\Sigma_{\is} \subseteq \Sigma_K$, we have also proven the `if' part of the second statement as a byproduct.

Next, to prove the `only if' part of the first statement, suppose that the data $(U_-,X)$ are informative for stabilization by state feedback. Let $K$ be such that $A+BK$ is stable for all $(A,B) \in \Sigma_{\is}$. By Lemma \ref{lemmaF0=0} we know that $$\im \begin{bmatrix}
I \\ K
\end{bmatrix} \subseteq 
\im \begin{bmatrix}
X_- \\ U_-
\end{bmatrix}.$$
This implies that $X_-$ has full row rank and there exists a right inverse $X_-^\dagger$ such that 
\begin{equation}
\label{eqK}
\begin{bmatrix}
I \\ K
\end{bmatrix} = \begin{bmatrix}
X_- \\ U_-
\end{bmatrix} X_-^\dagger.
\end{equation}
By \eqref{closedloop}, we obtain $A + BK = X_+ X_-^\dagger$, which shows that $X_+ X_-^\dagger$ is stable. This proves the `only if' part of the first statement. Finally, by \eqref{eqK}, the stabilizing feedback gain $K$ is indeed of the form $K = U_- X_-^\dagger$, which also proves the `only if' part of the second statement.
\EP 

Theorem \ref{t:algstab} gives a characterization of all data that are informative for stabilization by state feedback and provides a stabilizing controller. Nonetheless, the procedure to compute this controller might not be entirely satisfactory since it is not clear how to find a right inverse of $X_-$ that makes $X_+ X_-^\dagger$ stable. In general, $X_-$ has many right inverses, and $X_+ X_-^\dagger$ can be stable or unstable depending on the particular right inverse $X_-^\dagger$. To deal with this problem and to solve the design problem, we give a characterization of informativity for stabilization in terms of linear matrix inequalities (LMI's). The feasibility of such LMI's can be verified using standard methods. 

\begin{theorem}\label{t:lmistab}
	The data $(U_-,X)$ are informative for stabilization by state feedback if and only if there exists a matrix $\Theta \in \mathbb{R}^{T \times n}$ satisfying
	\begin{equation}
	\label{LMI/E}
	X_- \Theta = (X_- \Theta)^\top\quad\text{ and }\quad
	\begin{bmatrix}
	X_- \Theta & X_+ \Theta \\ \Theta^\top X_+^\top & X_- \Theta
	\end{bmatrix} > 0.
	\end{equation}
	Moreover, $K$ satisfies $\Sigma_{\is} \subseteq \Sigma_K$ if and only if $K = U_- \Theta (X_-\Theta)^{-1}$ for some matrix $\Theta$ satisfying \eqref{LMI/E}.
\end{theorem}
\begin{remark}
	To the best of our knowledge, LMI conditions for data-driven stabilization were first studied in \cite{DePersis2019}. In fact, the linear matrix inequality \eqref{LMI/E} is the same as that of \cite[Theorem 3]{DePersis2019}. However, an important difference is that the results in \cite{DePersis2019} assume that the input $ u $ is persistently exciting of sufficiently high order. In contrast, Theorem~\ref{t:lmistab}, as well as Theorem \ref{t:algstab}, do not require such conditions. The characterization \eqref{LMI/E} provides the minimal conditions on the data under which it is possible to obtain a stabilizing controller.
\end{remark}

\begin{example}
    Consider an unstable system of the form \eqref{e:disc}, where $A_s$ and $B_s$ are given by
	$$
	A_s = \begin{bmatrix}
	1.5 & 0 \\ 1 & 0.5
	\end{bmatrix}, \quad B_s = \begin{bmatrix}
	1 \\ 0
	\end{bmatrix}.
	$$
    We collect data from this system on a single time interval from $t = 0$ until $t = 2$, which results in the data matrices
    $$
    X = \begin{bmatrix}
    1 & 0.5 & -0.25 \\
    0 & 1 & 1
    \end{bmatrix}, \quad U_- = \begin{bmatrix}
    -1 & -1
    \end{bmatrix}.
    $$
    Clearly, the matrix $X_-$ is square and invertible, and it can be verified that
    $$
    X_+ X_-^{-1} = \begin{bmatrix}
    0.5 & -0.5 \\
    1 & 0.5
    \end{bmatrix}
    $$
    is stable, since its eigenvalues are $\half(1 \pm \sqrt{2}i)$. We conclude by Theorem \ref{t:algstab} that the data $(U_-,X)$ are informative for stabilization by state feedback. The same conclusion can be drawn from Theorem \ref{t:lmistab} since $$\Theta = \begin{bmatrix}
    1 & -1 \\ 0 & 2
    \end{bmatrix}$$ solves \eqref{LMI/E}. Next, we can conclude from either Theorem \ref{t:algstab} or Theorem \ref{t:lmistab} that the stabilizing feedback gain in this example is unique, and given by $K = U_-X_-^{-1} = \begin{bmatrix}
    -1 & -0.5
    \end{bmatrix}$. Finally, it is worth noting that the data are not informative for system identification. In fact, $(A,B) \in \Sigma_{\is}$ if and only if
    $$
    A = \begin{bmatrix}
    1.5+a_1 & 0.5a_1 \\ 1+ a_2 & 0.5+ 0.5a_2
    \end{bmatrix}, \quad B = \begin{bmatrix}
    1+a_1 \\ a_2
    \end{bmatrix}
    $$
    for some $a_1,a_2 \in \mathbb{R}$.
\end{example}

\noindent{\bf Proof of Theorem~\ref{t:lmistab}.}
To prove the `if' part of the first statement, suppose that there exists a $\Theta$ satisfying \eqref{LMI/E}. In particular, this implies that $X_- \Theta$ is symmetric positive definite. Therefore, $X_-$ has full row rank. By taking a Schur complement and multiplying by $-1$, we obtain 
$$
X_+ \Theta (X_- \Theta)^{-1} (X_- \Theta) (X_- \Theta)^{-1} \Theta^\top X_+^\top - X_- \Theta < 0.
$$ 
Since $X_- \Theta$ is positive definite, this implies that $X_+ \Theta (X_- \Theta)^{-1}$ is stable. In other words, there exists a right inverse $X_-^\dagger := \Theta (X_- \Theta)^{-1}$ of $X_-$ such that $X_+ X_-^\dagger$ is stable. By Theorem \ref{t:algstab}, we conclude that $(U_-,X)$ are informative for stabilization by state feedback, proving the `if' part of the first statement. Using Theorem \ref{t:algstab} once more, we see that $K := U_-\Theta (X_- \Theta)^{-1}$ stabilizes all systems in $\Sigma_{\is}$, which in turn proves the `if' part of the second statement.

Subsequently, to prove the `only if' part of the first statement, suppose that the data $(U_-,X)$ are informative for stabilization by state feedback. Let $K$ be any feedback gain such that $\Sigma_{\is} \subseteq \Sigma_K$. By Theorem \ref{t:algstab}, $X_-$ has full row rank and $K$ is of the form $K = U_- X_-^\dagger$, where $X_-^\dagger$ is a right inverse of $X_-$ such that $X_+ X_-^\dagger$ is stable. The stability of $X_+ X_-^\dagger$ implies the existence of a symmetric positive definite matrix $P$ such that 
\begin{align*}
(X_+ X_-^\dagger) P (X_+ X_-^\dagger)^\top - P < 0.
\end{align*}
Next, we define $\Theta := X_-^\dagger P$ and note that
$$
X_+ \Theta P^{-1} (X_+ \Theta)^\top - P < 0.
$$
Via the Schur complement we conclude that
$$
\begin{bmatrix}
P & X_+ \Theta \\ \Theta^\top X_+^\top & P
\end{bmatrix} > 0.
$$
Since $X_- X_-^\dagger = I$, we see that $P = X_- \Theta$, which proves the `only if' part of the first statement. Finally, by definition of $\Theta$, we have $X_-^\dagger = \Theta P^{-1} = \Theta (X_- \Theta)^{-1}$. Recall that $K = U_-X_-^\dagger$, which shows that $K$ is of the form $K = U_- \Theta (X_- \Theta)^{-1}$ for $\Theta$ satisfying \eqref{LMI/E}. This proves the `only if' part of the second statement and hence the proof is complete.
\EP

In addition to the stabilizing controllers discussed in Theorems \ref{t:algstab} and \ref{t:lmistab}, we may also look for a controller of the form $\bmu(t) = K \bmx(t)$ that stabilizes the system in \emph{finite time}. Such a controller is called a \emph{deadbeat controller} and is characterized by the property that $(A_s+B_sK)^t  x_0 = 0$ for all $t \geq n$ and all $x_0 \in \mathbb{R}^n$. Thus, $K$ is a deadbeat controller if and only if $A_s + B_sK$ is nilpotent. Now, for a given matrix $K$ define 
$$ 
\Sigma_K^{\text{nil}} := \{(A,B) \mid A+BK \text{ is nilpotent} \}.
$$
Then, analogous to the definition of informativity for stabilization by state feedback, we have the following definition of informativity for deadbeat control. 

\begin{definition}
	We say that the data $(U_-,X)$ are {\em informative for deadbeat control} if there exists a feedback gain $K$ such that $\Sigma_{\is} \subseteq \Sigma_K^{\text{nil}}$.
\end{definition}

Similarly to Theorem \ref{t:algstab}, we obtain the following necessary and sufficient conditions for informativity for deadbeat control.

\begin{theorem}
	\label{t:algdead}
	The data $(U_-,X)$ are informative for deadbeat control if and only if the matrix $X_-$ has full row rank and there exists a right inverse $X_-^\dagger$ of $X_-$ such that $X_+ X_-^\dagger$ is nilpotent.
	
	Moreover, if this condition is satisfied then the feedback gain $K := U_- X_-^\dagger$ yields a deadbeat controller, that is, $\Sigma_{\is} \subseteq \Sigma_K^{\text{nil}}$.
\end{theorem}

\begin{remark}
In order to compute a suitable right inverse $X_-^\dagger$ such that $X_+ X_-^\dagger$ is nilpotent, we can proceed as follows. Since $X_-$ has full row rank, we have $T\geq n$. We now distinguish two cases: $T=n$ and $T>n$. In the former case, $X_-$ is nonsingular and hence $X_+X_-\inv$ is nilpotent. In the latter case, there exist matrices $F\in\R^{T\times n}$ and $G\in\R^{T\times (T-n)}$ such that $\bbm F & G\ebm$ is nonsingular and $X_-\bbm F & G\ebm=\bbm I_n & 0_{n\times (T-n)}\ebm$. Note that $X_-^\dagger$ is a right inverse of $X_-$ if and only if $X_-^\dagger=F+GH$ for some $H\in\R^{(T-n)\times n}$. Finding a right inverse $X_-^\dagger$ such that $X_+ X_-^\dagger$ is nilpotent, therefore, amounts to finding $H$ such that $X_+F+X_+GH$ is nilpotent, i.e. has only zero eigenvalues. Such a matrix $H$ can be computed by invoking \cite[Thm. 3.29 and Thm. 3.32]{Trentelman2001} for the pair $(X_+F,X_+G)$ and the stability domain $\C_{\mathrm{g}}=\zset.$
\end{remark}

\subsection{Informativity for linear quadratic regulation}

Consider the discrete-time linear system \eqref{e:disc}. Let $x_{x_0,u}(\cdot)$ be the state sequence of \eqref{e:disc} resulting from the input $u(\cdot)$ and initial condition $x(0) = x_0$. We omit the subscript and simply write $x(\cdot)$ whenever the dependence on $x_0$ and $u$ is clear from the context. 

Associated to system \eqref{e:disc}, we define the quadratic cost functional
\beq\label{e:cost}
J(x_0,u)=\sum_{t=0}^\infty  x^\top(t) Q x(t) + u^\top(t) R u(t),
\eeq
where $Q = Q^\top$ is positive semidefinite and $R = R^\top$ is positive definite. Then, the linear quadratic regulator (LQR) problem is the following: 
\begin{problem}[LQR]
	Determine for every initial condition $x_0$ an input $u^*$, such that $\lim_{t\to\infty} x_{x_0,u^*}(t) = 0$, and the cost functional $J(x_0,u)$ is minimized under this constraint. 
\end{problem}
Such an input $u^*$ is called optimal for the given $x_0$. Of course, an optimal input does not necessarily exist for all $x_0$. We say that the linear quadratic regulator problem is {\em solvable\/} for $(A,B,Q,R)$ if for every $x_0$ there exists an input $u^*$ such that
\begin{enumerate}
	\item The cost $J(x_0,u^*)$ is finite.
	\item The limit $\lim_{t\to\infty}x_{x_0,u^*}(t)=0$. 
	\item The input $u^*$ minimizes the cost functional, i.e., 
	\[J(x_0,u^*)\leq J(x_0,\baru)\]
	for all $\baru$ such that $\lim_{t\to\infty}x_{x_0,\baru}(t)=0$.
\end{enumerate}

In the sequel, we will require the notion of observable eigenvalues. Recall from e.g. \cite[Section 3.5]{Trentelman2001} that an eigenvalue $\lambda$ of $A$ is $(Q,A)$-observable if 
\[ \rank \begin{pmatrix} A-\lambda I \\ Q \end{pmatrix}=n.\] 

The following theorem provides necessary and sufficient conditions for the solvability of the linear quadratic regulator problem for $(A,B,Q,R)$. This theorem is the discrete-time analogue to the continuous-time case stated in \cite[Theorem 10.18]{Trentelman2001}. 

\begin{theorem}\label{t:Harry}
	Let $Q = Q^\top$ be positive semidefinite and $R = R^\top$ be positive definite. Then the following statements hold:	
	\begin{enumerate}[(i)]
		\item If $(A,B)$ is stabilizable, there exists a unique largest real symmetric solution $P^+$ to the discrete-time algebraic Riccati equation (DARE) 
		\begin{equation}
		\label{dare}
		P = A^\top PA-A^\top PB(R+B^\top P B)\inv B^\top  P A+Q,
		\end{equation}
		in the sense that $P^+ \geq P$ for every real symmetric $P$ satisfying \eqref{dare}. The matrix $P^+$ is positive semidefinite.
		\item If, in addition to stabilizability of $(A,B)$, every eigenvalue of $A$ on the unit circle is $(Q,A)$-observable then for every $x_0$ a unique optimal input $u^*$ exists. Furthermore, this input sequence is generated by the feedback law $\bmu = K \bmx$, where
		\begin{equation}
		\label{optgain}
		K := -(R+B^\top P^+ B)\inv B^\top  P^+ A.
		\end{equation}
		Moreover, the matrix $A+BK$ is stable. 
		\item In fact, the linear quadratic regulator problem is solvable for $(A,B,Q,R)$ if and only if $(A,B)$ is stabilizable and every eigenvalue of $A$ on the unit circle is $(Q,A)$-observable. 
	\end{enumerate}
\end{theorem}

If the LQR problem is solvable for $(A,B,Q,R)$, we say that $K$ given by \eqref{optgain} is the optimal feedback gain for $(A,B,Q,R)$. 

Now, for any given $K$ we define $\Sigma^{Q,R}_{K}$ as the set of all systems of the form \eqref{e:disc} for which $K$ is the optimal feedback gain corresponding to $Q$ and $R$, that is,
\[ \Sigma_K^{Q,R}:=\set{(A,B)}{K \text{ is the optimal gain for }(A,B,Q,R)}.\]
This gives rise to another notion of informativity in line with Definition~\ref{def:par informativity}. Again, let $\Sigma_{\is}$ be given by \eqref{eq:sigma is}.

\begin{definition}
	Given matrices $Q$ and $R$, we say that the data $(U_-,X)$ are \emph{informative for linear quadratic regulation} if there exists $K$ such that $\Sigma_{\is} \subseteq \Sigma^{Q,R}_{K}$.
\end{definition}

In order to provide necessary and sufficient conditions for the corresponding informativity problem, we need the following auxiliary lemma.

\begin{lemma}\label{l: same P works for all}
Let $Q=Q^\top$ be positive semidefinite and $R=R^\top$ be positive definite. Suppose the data $(U_-,X)$ are informative for linear quadratic regulation. Let $K$ be such that $\Sigma_{\is} \subseteq \Sigma^{Q,R}_{K}$. Then, there exist a square matrix $M$ and a symmetric positive semidefinite matrix $P^+$ such that for all $(A,B)\in\Sigma_{\is}$
\begin{align}
&\!\!\!\!\!M=A+BK,\label{e: define M}\\
&\!\!\!\!\!P^+\!= A^\top\! P^+\!A\! -\! A^\top\! P^+\!B(R + B^\top\! P^+\! B)\inv B^\top\!  P^+\! A + Q,\!\!\label{e:dare aux}\\
&\!\!\!\!\!P^+-M^\top P^+M=K^\top RK+Q,\label{e:lyap aux}\\
&\!\!\!\!\!K=-(R+B^\top P^+ B)\inv B^\top  P^+ A.\label{e:uni K aux}
\end{align}
\end{lemma}

\BP
%Since the data $(U_-,X)$ are informative for linear quadratic regulation, there exists a $K$ that is the optimal feedback gain for $(A,B,Q,R)$ for all $(A,B)\in\Sigma_{\is}$. It is noteworthy that $K$ is one and the same for all $(A,B)\in\Sigma_{\is}$. 
Since the data $(U_-,X)$ are informative for linear quadratic regulation, $A+BK$ is stable for every $(A,B)\in\Sigma_{\is}$. By Lemma~\ref{lemmaF0=0}, this implies that $A_0 + B_0 K = 0$ for all $(A_0,B_0) \in \Sigma_{\is}^0$. Thus, there exists $M$ such that $M=A+BK$ for all $(A,B) \in \Sigma_{\is}$. For the rest, note that Theorem~\ref{t:Harry} implies that for every $(A,B)\in\Sigma_{\is}$ there exists $P^+_{(A,B)}$ satisfying the DARE 
\begin{equation}\label{e:dare for p plus}
\begin{split} 
&P^+_{(A,B)}= A^\top P^+_{(A,B)}A  \\ &-A^\top P^+_{(A,B)}B(R+B^\top P^+_{(A,B)} B)\inv B^\top  P^+_{(A,B)} A+Q 
\end{split}
\end{equation}
such that 
\begin{equation}\label{e: exp K}
K = -(R+B^\top P^+_{(A,B)} B)\inv B^\top  P^+_{(A,B)} A.
\end{equation}
It is important to note that, although $K$ is independent of the choice of $(A,B)$, the matrix $P^+_{(A,B)}$ might depend on $(A,B)$. We will, however, show that also $P^+_{(A,B)}$ is independent of the choice of $(A,B)$. 

By rewriting \eqref{e:dare for p plus}, we see that
\begin{equation}\label{e:dare alt}
P^+_{(A,B)}-M^\top P^+_{(A,B)}M=K^\top RK+Q.
\end{equation}
Since $M$ is stable, $P^+_{(A,B)}$ is the unique solution to the discrete-time Lyapunov equation \eqref{e:dare alt}, see e.g. \cite[Section 6]{Simoncini2016}. Moreover, since $M$ and $K$ do not depend on the choice of $(A,B)\in\Sigma_{\is}$, it indeed follows that $P^+_{(A,B)}$ does not depend on $(A,B)$. 
It follows from \eqref{e:dare for p plus}--\eqref{e:dare alt} that $P^+:=P^+_{(A,B)}$ satisfies \eqref{e:dare aux}--\eqref{e:uni K aux}.
\EP

The following theorem solves the informativity problem for linear quadratic regulation.

\begin{theorem}\label{t:LQinform}
	Let $Q=Q^\top$ be positive semidefinite and $R=R^\top$ be positive definite. Then, the data $(U_-,X)$ are informative for linear quadratic regulation if and only if at least one of the following two conditions hold:
	\begin{enumerate}[(i)]
		\item\label{cond:a} The data $(U_-,X)$ are informative for system identification, that is, $\Sigma_{\is}=\pset{(A_s,B_s)}$, and the linear quadratic regulator problem is solvable for $(A_s,B_s,Q,R)$. In this case, the optimal feedback gain $K$ is of the form \eqref{optgain} where $P^+$ is the largest real symmetric solution to \eqref{dare}.
		\item\label{cond:b} For all $(A,B) \in \Sigma_{\is}$ we have $A=A_s$. Moreover, $A_s$ is stable, $QA_s = 0$, and the optimal feedback gain is given by $K = 0$. 
	\end{enumerate}
\end{theorem}
\begin{remark} Condition \eqref{cond:b} of Theorem \ref{t:LQinform} is a pathological case in which $A$ is stable and $QA = 0$ for all matrices $A$ that are compatible with the data. Since $x(t) \in \im A$ for all $t > 0$, we have $Q x(t) = 0$ for all $t > 0$ if the input function is chosen as $u = 0$. Additionally, since $A$ is stable, this shows that the optimal input is equal to $u^* = 0$. If we set aside condition \eqref{cond:b}, the implication of Theorem \ref{t:LQinform} is the following: if the data are informative for linear quadratic regulation they are also informative for system identification.

At first sight, this might seem like a negative result in the sense that data-driven LQR is only possible with data that are also informative enough to uniquely identify the system. However, at the same time, Theorem \ref{t:LQinform} can be viewed as a positive result in the sense that it provides fundamental justification for the data conditions imposed in e.g. \cite{DePersis2019}. Indeed, in \cite{DePersis2019} the data-driven infinite horizon LQR problem\footnote{Note that the authors of \cite{DePersis2019} formulate this problem as the minimization of the $H_2$-norm of a certain transfer matrix.} is solved using input/state data under the assumption that the input is persistently exciting of sufficiently high order. Under the latter assumption, the input/state data are informative for system identification, i.e., the matrices $A_s$ and $B_s$ can be uniquely determined from data. Theorem \ref{t:LQinform} justifies such a strong assumption on the richness of data in data-driven linear quadratic regulation.

The data-driven \emph{finite} horizon LQR problem was solved under a persistency of excitation assumption in \cite{Markovsky2007}. Our results suggest that also in this case informativity for system identification is necessary for data-driven LQR, although further analysis is required to prove this claim.
\end{remark}

\noindent{\bf Proof of Theorem~\ref{t:LQinform}.}
We first prove the `if' part. Sufficiency of the condition \eqref{cond:a} readily follows from Theorem~\ref{t:Harry}. To prove the sufficiency of the condition \eqref{cond:b}, assume that the matrix $A$ is stable and $QA = 0$ for all $(A,B) \in \Sigma_{\is}$. By the discussion following Theorem \ref{t:LQinform}, this implies that $u^* = 0$ for all $(A,B) \in \Sigma_{\is}$. Hence, for $K = 0$ we have $\Sigma_{\is} \subseteq \Sigma_K^{Q,R}$, i.e., the data are informative for linear quadratic regulation.

To prove the `only if' part, suppose that the data $(U_-,X)$ are informative for linear quadratic regulation. From Lemma~\ref{l: same P works for all}, we know that there exist $M$ and $P^+$ satisfying \eqref{e: define M}--\eqref{e:uni K aux} for all $(A,B)\in\Sigma_{\is}$. By substituting \eqref{e:uni K aux} into \eqref{e:dare aux} and using \eqref{e: define M}, we obtain 
\beq\label{e:A0 part}
A^\top P^+M=P^+-Q.
\eeq
In addition, it follows from \eqref{e:uni K aux} that $-(R+B^\top P^+B)K=B^\top P^+A$. By using \eqref{e: define M}, we have
\beq\label{e:B0 part}
B^\top P^+M=-RK.
\eeq
Since \eqref{e:A0 part} and \eqref{e:B0 part} hold for all $(A,B)\in\Sigma_{\is}$, we have that
$$
\bbm A_0^\top\\B^\top_0\ebm P^+M=0
$$
for all $(A_0,B_0)\in\Sigma_{\is}^0$. Note that $(FA_0,FB_0)\in\Sigma_{\is}^0$ for all $F\in\R^{n\times n}$ whenever $(A_0,B_0)\in\Sigma_{\is}^0$. This means that
$$
\bbm A_0^\top\\B^\top_0\ebm F^\top P^+M=0
$$
for all $F\in\R^{n\times n}$. Therefore, either $\bbm A_0&B_0\ebm=0$ for all $(A_0,B_0)\in\Sigma_{\is}^0$ or $P^+M=0$. The former is equivalent to $\Sigma_{\is}^0=\zset$. In this case, we see that the data $(U_-,X)$ are informative for system identification, equivalently $\Sigma_{\is}=\pset{(A_s,B_s)}$, and the LQR problem is solvable for $(A_s,B_s,Q,R)$. Therefore, condition \eqref{cond:a} holds. On the other hand, if 
$P^+M=0$ then we have
\begin{align*}
0=P^+M &= P^+(A+BK)\\
&= P^+\big(A-B(R+B^\top P^+ B)\inv B^\top  P^+ A\big) \\
&= \big(I-P^+B(R+B^\top P^+B)\inv B^\top\big)P^+A.
\end{align*}
for all $(A,B)\in\Sigma_{\is}$. From the identity
$$
(I+P^+B R\inv B^\top)\inv=I-P^+B(R+B^\top P^+B)\inv B^\top,
$$
we see that $P^+A=0$ for all $(A,B)\in\Sigma_{\is}$. Then, it follows from \eqref{e:uni K aux} that $K = 0$. Since $A_0+B_0K=0$ for all $(A_0,B_0)\in\Sigma_{\is}^0$ due to Lemma~\ref{lemmaF0=0}, we see that $A_0$ must be zero. Hence, 
we have $A=A_s$ for all $(A,B) \in \Sigma_{\is}$ and $A_s$ is stable. Moreover, it follows from \eqref{e:A0 part} that $P^+ = Q$. Therefore, $QA_s = 0$. In other words, condition \eqref{cond:b} is satisfied, which proves the theorem. \EP

Theorem \ref{t:LQinform} gives necessary and sufficient conditions under which the data are informative for linear quadratic regulation. However, it might not be directly clear how these conditions can be verified given input/state data. Therefore, in what follows we rephrase the conditions of Theorem \ref{t:LQinform} in terms of the data matrices $X$ and $U_-$.

\begin{theorem}
	\label{t:LQinform2}
	Let $Q=Q^\top$ be positive semidefinite and $R=R^\top$ be positive definite. Then, the data $(U_-,X)$ are informative for linear quadratic regulation if and only if at least one of the following two conditions hold:
	\begin{enumerate}[(i)]
		\item\label{cond:a2} The data $(U_-,X)$ are informative for system identification. Equivalently, there exists $\begin{bmatrix} V_1 & V_2\end{bmatrix}$ such that \eqref{eq:V1 V2} holds. Moreover, the linear quadratic regulator problem is solvable for $(A_s,B_s,Q,R)$, where $A_s=X_+V_1$ and $B_s=X_+V_2$.
		\item\label{cond:b2} There exists $\Theta \in \mathbb{R}^{T \times n}$ such that $	X_- \Theta = (X_- \Theta)^\top$, 	$U_- \Theta = 0 $, 
		\begin{equation}
		\label{e:LMI/E/K/Q}
		\begin{bmatrix}
		X_- \Theta & X_+ \Theta \\ \Theta^\top X_+^\top & X_- \Theta
		\end{bmatrix} > 0.
		\end{equation}
		and $ QX_+\Theta = 0$. 
	\end{enumerate}
\end{theorem}

\BP
The equivalence of condition \eqref{cond:a} of Theorem \ref{t:LQinform} and condition \eqref{cond:a2} of Theorem \ref{t:LQinform2} is obvious. It remains to be shown that condition \eqref{cond:b} of Theorem \ref{t:LQinform} and condition \eqref{cond:b2} of Theorem \ref{t:LQinform2} are equivalent as well. To this end, suppose that there exists a matrix $\Theta \in \mathbb{R}^{T \times n}$ such that the conditions of \eqref{cond:b2} holds. By Theorem \ref{t:lmistab}, we have $\Sigma_{\is} \subseteq \Sigma_K$ for $K = 0$, that is, $A$ is stable for all $(A,B) \in \Sigma_{\is}$. In addition, note that
\begin{equation}
\label{eqQ}
Q X_+ \Theta (X_-\Theta)^{-1} = Q \begin{bmatrix}
A & B
\end{bmatrix} \begin{bmatrix}
X_- \\ U_-
\end{bmatrix} \Theta (X_-\Theta)^{-1} = QA
\end{equation}  
for all $(A,B) \in \Sigma_{\is}$. This shows that $QA=0$ and therefore that condition \eqref{cond:b} of Theorem \ref{t:LQinform} holds. Conversely, suppose that $A$ is stable and $QA = 0$ for all $(A,B) \in \Sigma_{\is}$. This implies that $K = 0$ is a stabilizing controller for all $(A,B) \in \Sigma_{\is}$. By Theorem \ref{t:lmistab}, there exists a matrix $\Theta \in \mathbb{R}^{T \times n}$ satisfying the first three conditions of \eqref{cond:b2}. Finally, it follows from $QA = 0$ and \eqref{eqQ} that $\Theta$ also satisfies the fourth equation of \eqref{cond:b2}. This proves the theorem. 
\EP 

\subsection{From data to LQ gain}
In this section our goal is to devise a method in order to compute the optimal feedback gain $K$ directly from the data. For this, we will employ ideas from the study of Riccati inequalities (see e.g \cite{Ran1988}).
%
% The main idea of our approach will be to replace the Riccati inequality by the linear matrix inequality $\mathcal{L}(P) \leq 0$, where
%\begin{equation}
%\label{dataineq}
%\mathcal{L}(P) := \xmt P\xm-\xpt P\xp -\xmt Q\xm-\umt R\um.
%\end{equation}
%Note that the linear operator $\mathcal{L}$ is completely defined by the data matrices $X$ and $U_-$, and the weight matrices $Q$ and $R$. 

The following theorem asserts that $P^+$ as in Lemma~\ref{l: same P works for all} can be found as the unique solution to an optimization problem involving only the data. Furthermore, the optimal feedback gain $K$ can subsequently be found by solving a set of linear equations. 

\begin{theorem}
	\label{t:LQgaindata}
Let $Q = Q^\top \geq 0$ and $R = R^\top > 0$. Suppose that the data $(U_-,X)$ are informative for linear quadratic regulation. Consider the linear operator $P\mapsto\calL(P)$ defined by
$$
\mathcal{L}(P) := \xmt P\xm-\xpt P\xp -\xmt Q\xm-\umt R\um.
$$
Let $P^+$ be as in Lemma~\ref{l: same P works for all}. The following statements hold:
	\begin{enumerate}[(i)]
		\item \label{semidefiniteprogram} The matrix $P^+$ is equal to the unique solution to the optimization problem
		\begin{align*}
		\text{ maximize } \: &\trace P \\
		\text{ subject to } \: &P = P^\top \geq 0 
		\,\,\text{ and }\,\,  \mathcal{L}(P) \leq 0.
		\end{align*} 
		\item There exists a right inverse $X_-^\dagger$ of $X_-$ such that
			\begin{align}
			\label{eq1}
			\mathcal{L}(P^+) X_-^\dagger &= 0.
			\end{align}
		Moreover, if $X_-^\dagger$ satisfies \eqref{eq1}, then the optimal feedback gain is given by $K = U_- X_-^\dagger$.
	\end{enumerate}
	
\end{theorem}

\begin{remark}
	From a design viewpoint, the optimal feedback gain $K$ can be found in the following way. First solve the semidefinite program in Theorem \ref{t:LQgaindata}\eqref{semidefiniteprogram}. Subsequently, compute a solution $X_-^\dagger$ to the linear equations $X_- X_-^\dagger = I$ and \eqref{eq1}. Then, the optimal feedback gain is given by $K = U_- X_-^\dagger$.
\end{remark}

\begin{remark}
	
	The data-driven LQR problem was first solved using semidefinite programming in \cite[Theorem 4]{DePersis2019}. There, the optimal feedback gain was found by minimizing the trace of a weighted sum of two matrix variables, subject to two LMI constraints. The semidefinite program in Theorem \ref{t:LQgaindata} is attractive since the dimension of the unknown $P$ is (only) $n \times n$. In comparison, the dimensions of the two unknowns in \cite[Theorem 4]{DePersis2019} are $T \times n$ and $m \times m$, respectively. In general, the number of samples $T$ is much larger\footnote{In fact, this is always the case under the persistency of excitation conditions imposed in \cite{DePersis2019} as such conditions can only be satisfied provided that $T \geq nm + n + m$.} than $n$. An additional attractive feature of Theorem \ref{t:LQgaindata} is that $P^+$ is obtained from the data. This is useful since the minimal cost associated to any initial condition $x_0$ can be computed as $x_0^\top P^+ x_0$.
	
	The data-driven LQR approach in \cite{daSilva2019} is quite different from Theorem \ref{t:LQgaindata} since the solution to the Riccati equation is approximated using a batch-form solution to the \emph{Riccati difference equation}. A similar approach was used in \cite{Skelton1994,Furuta1995,Shi1998,Aangenent2005} for the \emph{finite horizon} data-driven LQR/LQG problem. In the setup of \cite{daSilva2019}, the approximate solution to the Riccati equation is exact only if the number of data points tends to infinity. The main difference between our approach and the one in \cite{daSilva2019} is hence that the solution $P^+$ to the Riccati equation can be obtained exactly from \emph{finite} data via Theorem \ref{t:LQgaindata}. 
	
\end{remark}

\noindent{\bf Proof of Theorem~\ref{t:LQgaindata}.}
We begin with proving the first statement. Note that 
$$
\mathcal{L}(P)=
\bbm\xm\\\um\ebm^\top\!\!\bbm P-A^\top PA-Q & -A^\top PB\\
-B^\top PA & -(R+B^\top PB)\ebm\bbm\xm\\\um\ebm
$$
for all $(A,B)\in\Sigma_{\is}$. We claim that the following implication holds:
\beq
P = P^\top \geq 0 \,\,\text{ and }\,\,  \mathcal{L}(P) \leq 0\implies P^+\geq P.\label{e:claim p to lp}
\eeq
To prove this claim, let $P$ be such that $P= P^\top \geq 0$ and $\calL(P)\leq0$. Since the data are informative for linear quadratic regulation, they are also informative for stabilization by state feedback. Therefore, the optimal feedback gain $K$ satisfies
$$
\im\bbm I \\ K\ebm\subseteq\im\bbm X_-\\U_-\ebm
$$
due to Lemma~\ref{lemmaF0=0}. Therefore, the above expression for $\calL(P)$ implies that
$$
\bbm I \\ K\ebm^\top\!\!\bbm P-A^\top PA-Q & -A^\top PB\\
-B^\top PA & -(R+B^\top PB)\ebm\bbm I \\ K\ebm\leq 0
$$
for all $(A,B)\in\Sigma_{\is}$. This yields
$$
P - M^\top P M\leq K^\top RK+ Q
$$
where $M$ is as in Lemma~\ref{l: same P works for all}. By subtracting this from \eqref{e:lyap aux}, we obtain
$$
(P^+-P) - M^\top (P^+-P) M\geq 0.
$$
Since $M$ is stable, this discrete-time Lyapunov inequality implies that $P^+-P\geq 0$ and hence $P^+\geq P$. This proves the claim \eqref{e:claim p to lp}.

Note that $R+B^\top P^+B$ is positive definite. Then, it follows from \eqref{e:dare aux} that
$$
\bbm P^+-A^\top P^+A-Q & -A^\top P^+B\\
-B^\top P^+A & -(R+B^\top P^+B)\ebm\leq 0
$$
via a Schur complement argument. Therefore, $\calL(P^+)\leq 0$. Since $P^+\geq P$, we have $\trace P^+\geq \trace P$. Together with \eqref{e:claim p to lp}, this shows that $P^+$ is a solution to the optimization problem stated in the theorem.

Next, we prove uniqueness. Let $\barP$ be another solution of the optimization problem. Then,  we have that $\barP=\barP^\top\geq 0$, $\calL(\barP)\leq 0$, and $\trace \barP=\trace P^+$. From \eqref{e:claim p to lp}, we see that $P^+ \geq \barP$. In particular, this implies that $(P^+)_{ii} \geq \barP_{ii}$ for all $i$. Together with $\trace \barP=\trace P^+$, this implies that $(P^+)_{ii} = \barP_{ii}$ for all $i$. Now, for any $i$ and $j$, we have
\begin{align*}
(e_i - e_j)^\top P^+ (e_i - e_j) &\geq (e_i - e_j)^\top \barP (e_i - e_j) \text{ and } \\
(e_i + e_j)^\top P^+ (e_i + e_j) &\geq (e_i + e_j)^\top \barP (e_i + e_j),
\end{align*}
where $e_i$ denotes the $i$-th standard basis vector. This leads to $(P^+)_{ij} \leq \barP_{ij}$ and $(P^+)_{ij} \geq \barP_{ij}$, respectively. We conclude that $(P^+)_{ij} = \barP_{ij}$ for all $i,j$. This proves uniqueness. 

Finally, we prove the second statement. It follows from \eqref{e:dare aux} and \eqref{e:uni K aux} that
\begin{equation}
\label{LP+}
\mathcal{L}(P^+) = - \left(U_- - K X_-\right)^\top (R+B^\top P^+B)
\left(U_- - K X_-\right).
\end{equation}
The optimal feedback $K$ is stabilizing, therefore it follows from Theorem \ref{t:algstab} that $K$ can be written as $K = U_- \Gamma$, where $\Gamma$ is some right inverse of $X_-$. Note that this implies the existence of a right inverse $X_-^\dagger$ of $X_-$ satisfying \eqref{eq1}. Indeed, $X_-^\dagger := \Gamma$ is such a matrix by \eqref{LP+}. Moreover, if $X_-^\dagger$ is a right inverse of $X_-$ satisfying \eqref{eq1} then $(U_- - KX_-)X_-^\dagger = 0$ by \eqref{LP+} and positive definiteness of $R$. We conclude that the optimal feedback gain is equal to $K = U_- X_-^\dagger$, which proves the second statement. 
\EP

\section{Control using input and output data}\label{sec:output}
In this section, we will consider problems where the output \textit{does} play a role. In particular, we will consider the problem of stabilization by dynamic measurement feedback. We will first consider this problem based on input, state and output measurements. Subsequently, we will turn our attention to the case of input/output data. 

Consider the `true' system 
\begin{subequations}\label{inoutsys}
	\begin{align}
	\label{inoutsys1}
	\bmx(t+1) &= A_s\bmx(t) + B_s\bmu(t) \\
	\label{inoutsys2}
	\bmy(t) &= C_s\bmx(t) + D_s\bmu(t).
	\end{align}
\end{subequations}
We want to design a stabilizing dynamic controller of the form%
\begin{subequations}\label{compensator}
	\begin{align}
	\label{compensator1}
	\bmw(t+1) &= K \bmw(t) + L\bmy(t) \\
	\label{compensator2}
	\bmu(t) &= M\bmw(t)
	\end{align}
\end{subequations}
such that the closed-loop system, given by 
\begin{equation*}
\begin{bmatrix}
\bmx(t+1) \\ \bmw(t+1)
\end{bmatrix} = \begin{bmatrix}
A_s & B_sM \\ LC_s & K + LD_sM
\end{bmatrix} \begin{bmatrix}
\bmx(t) \\ \bmw(t)
\end{bmatrix},
\end{equation*}
is stable. This is equivalent to the condition that
\begin{equation}
\label{closedloopKLM}
\begin{bmatrix}
A_s & B_sM \\ LC_s & K + LD_sM
\end{bmatrix}
\end{equation}
is a stable matrix. 
\subsection{Stabilization using input, state and output data}\label{sec:iso measurements}
Suppose that we collect input/state/output data on $\ell$ time intervals $\pset{0,1,\ldots,T_i}$ for $i=1,2,\ldots,q$. Let $U_-,X,X_-,$ and $X_+$ be defined as in \eqref{eq: UXdatanew} and let $Y_-$ be defined in a similar way as $U_-$.
Then, we have
\begin{equation}
\label{dataeq}
\begin{bmatrix}
X_+ \\ Y_-
\end{bmatrix} = \begin{bmatrix}
A_s & B_s \\ C_s & D_s
\end{bmatrix} \begin{bmatrix}
X_- \\ U_-
\end{bmatrix}
\end{equation}
relating the data and the `true' system \eqref{inoutsys}. The set of all systems that are consistent with these data is then given by: 
\begin{equation}\label{eq:sigma iso}
\Sigma_{\iso} := \left\{ (A,B,C,D) \mid \begin{bmatrix} X_+ \\ Y_- \end{bmatrix} = \begin{bmatrix} A & B \\ C & D \end{bmatrix} \begin{bmatrix} X_- \\ U_- \end{bmatrix} \right\}.
\end{equation}
In addition, for given $K$, $L$ and $M$, we define the set of systems that are stabilized by the dynamic controller \eqref{compensator} by
\begin{equation*}
\Sigma_{K,L,M} := \left\{ (A,B,C,D) \mid \begin{bmatrix} A & BM \\ LC & K + LDM \end{bmatrix} \text{ is stable} \right\}.
\end{equation*}

Subsequently, in line with Definition~\ref{def:par informativity}, we consider the following notion of informativity:

\begin{definition}
	We say the data $(U_-,X,Y_-)$ are \emph{informative for stabilization by dynamic measurement feedback} if there exist matrices $K$, $L$ and $M$ such that $\Sigma_{\iso} \subseteq \Sigma_{K,L,M}$.
\end{definition}

As in the general case of informativity for control, we consider two consequent problems: First, to characterize informativity for stabilization in terms of necessary and sufficient conditions on the data and next to design a controller based on these data. To aid in solving these problems, we will first investigate the case where $U_-$ does not have full row rank. In this case, we will show that the problem can be `reduced' to the full row rank case.

For this, we start with the observation that any $U_-\in \mathbb{R}^{m\times T}$ of row rank $k < m$ can be decomposed as $U_-= S\hat{U}_- $, where $S$ has full column rank and $\hat{U}_-\in\mathbb{R}^{k\times T}$ has full row rank. We now have the following lemma:

\begin{lemma}
	\label{l:equiv}
	Consider the data $(U_-,X,Y_-)$ and the corresponding set $\Sigma_{\iso}$. Let $S$ be a matrix of full column rank such that $U_- = S \hat{U}_-$ with $\hat{U}_-$ a matrix of full row rank. Let $S^\dagger$ be a left inverse of $S$.
	
	Then the data $(U_-,X,Y_-)$ are informative for stabilization by dynamic measurement feedback if and only if the data $(\hat{U}_-,X,Y_-)$ are informative for stabilization by dynamic measurement feedback. 
	
	In particular, if we let $\hat{\Sigma}_{\iso}$ be the set of systems consistent with the `reduced' data set $(\hat{U}_-,X,Y_-)$, and if $\hat{K}$ $\hat{L}$ and $\hat{M}$ are real matrices of appropriate dimensions, then: 
		\begin{align}
		\Sigma_{\iso} \subseteq \Sigma_{K,L,M} &\implies \hat{\Sigma}_{\iso} \subseteq \Sigma_{K,L,S^\dagger M}, \label{eq: iso in klm1} \\
		\hat{\Sigma}_{\iso} \subseteq \Sigma_{\hat{K},\hat{L},\hat{M}} &\implies \Sigma_{\iso} \subseteq \Sigma_{\hat{K},\hat{L},S\hat{M}}. \label{eq: iso in klm2} 
		\end{align}	
\end{lemma}

\BP
	First note that
	\[\hat{\Sigma}_{\iso} = \left\{ (\hat{A},\hat{B},\hat{C},\hat{D}) \mid \begin{bmatrix}
	X_+ \\ Y_-
	\end{bmatrix} = \begin{bmatrix}
	\hat{A} & \hat{B} \\ \hat{C} & \hat{D}
	\end{bmatrix} \begin{bmatrix}
	X_- \\ \hat{U}_-
	\end{bmatrix} \right\}. \]	
	We will start by proving the following two implications: 
		\begin{align}
		\!(A,B,C,D) \in \Sigma_{\iso} &\!\!\!\implies\!\! (A,BS,C,DS) \in \hat{\Sigma}_{\iso}, \label{imp1} \\ \!(\hat{A},\hat{B},\hat{C},\hat{D}) \in \hat{\Sigma}_{\iso} &\!\!\!\implies\!\! (\hat{A},\hat{B}S^\dagger,\hat{C},\hat{D}S^\dagger) \in \Sigma_{\iso}. \label{imp2}
		\end{align}
	To prove implication \eqref{imp1}, assume that $(A,B,C,D) \in \Sigma_{\iso}$. Then, by definition
	\[\begin{bmatrix} 	X_+ \\ Y_-\end{bmatrix} = \begin{bmatrix}	A & B \\ C & D\end{bmatrix} \begin{bmatrix} X_- \\ U_-\end{bmatrix}.\]
	From the definition of $S$, we have $U_- = S\hat{U}_-$. Substitution of this results in
	\[\begin{bmatrix} 	X_+ \\ Y_-\end{bmatrix} = \begin{bmatrix}	A & B \\ C & D\end{bmatrix} \begin{bmatrix} X_- \\ S \hat{U}_-\end{bmatrix}=\begin{bmatrix}	A & BS \\ C & DS\end{bmatrix} \begin{bmatrix} X_- \\ \hat{U}_-\end{bmatrix}.\]
	This implies that $(A,BS,C,DS) \in \hat{\Sigma}_{\iso}$. The implication \eqref{imp2} can be proven similarly by substitution of $\hat{U}_- = S^\dagger U_-$. 
	
	To prove the lemma, suppose that the data $(U_-,X,Y_-)$ are informative for stabilization by dynamic measurement feedback. This means that there exist $K$, $L$, and $M$ such that $$\begin{bmatrix}
	A & BM \\ LC & K + LDM
	\end{bmatrix}$$ is stable for all $(A,B,C,D) \in \Sigma_{\iso}$. In particular, if $(\hat{A},\hat{B},\hat{C},\hat{D}) \in \hat{\Sigma}_{\iso}$ then  $(\hat{A},\hat{B}S^\dagger,\hat{C},\hat{D}S^\dagger ) \in \Sigma_{\iso}$ by \eqref{imp2}. This means that the matrix 
	\[\begin{bmatrix} \hat{A} & \hat{B}S^\dagger M \\ L\hat{C} & K + L\hat{D}S^\dagger M	\end{bmatrix}\]
	is stable for all $(\hat{A},\hat{B},\hat{C},\hat{D}) \in \hat{\Sigma}_{\iso}$. In other words, $\hat{\Sigma}_{\iso} \subseteq \Sigma_{K,L,S^\dagger M}$ and hence implication \eqref{eq: iso in klm1} holds and the data $(\hat{U}_-,X,Y_-)$ are informative for stabilization by dynamic measurement feedback. The proofs of \eqref{eq: iso in klm2} and the `if' part of the theorem are analogous and hence omitted.
\EP

We will now solve the informativity and design problems under the condition that $U_-$ has full row rank. 

\begin{theorem}
	\label{t:stabUXY}
	Consider the data $(U_-,X,Y_-)$ and assume that $U_-$ has full row rank. Then $(U_-,X,Y_-)$ are informative for stabilization by dynamic measurement feedback if and only if the following  conditions are satisfied:
	\begin{enumerate}[(i)]
		\item We have
		\[ \rank\begin{bmatrix} X_- \\ U_- \end{bmatrix} = n+m.\]
		Equivalently, there exists $\begin{bmatrix} V_1 & V_2\end{bmatrix}$ such that \eqref{eq:V1 V2} holds. This means that $$\Sigma_{\iso} = \{ (X_+V_1, X_+V_2,Y_-V_1,Y_-V_2)\}.$$  \label{1}
		\item The pair $(X_+ V_1, X_+ V_2)$ is stabilizable and $(Y_- V_1, X_+ V_1)$ is detectable. \label{2}
	\end{enumerate}
	Moreover, if the above conditions are satisfied, a stabilizing controller $(K,L,M)$ can be constructed as follows:
	\begin{enumerate}[(a)]
		\item Select a matrix $M$ such that $X_+ (V_1 + V_2 M)$ is stable. \label{a}
		\item Choose a matrix $L$ such that $(X_+ - L Y_-) V_1$ is stable. \label{b}
		\item Define $K := (X_+ - LY_-)(V_1 + V_2 M)$. \label{c}
	\end{enumerate} 
\end{theorem}
\begin{remark}
	Under the condition that $U_-$ has full row rank, Theorem \ref{t:stabUXY} asserts that in order to construct a stabilizing dynamic controller, it is necessary that the data are rich enough to identify the system matrices $A_s,B_s,C_s$ and $D_s$ uniquely. The controller proposed in \eqref{a}, \eqref{b}, \eqref{c} is a so-called \emph{observer-based} controller, see e.g. \cite[Section 3.12]{Trentelman2001}. The feedback gains $M$ and $L$ can be computed using standard methods, for example via pole placement or LMI's.   
\end{remark}

\noindent{\bf Proof of Theorem~\ref{t:stabUXY}.}	To prove the `if' part, suppose that conditions \eqref{1} and \eqref{2} are satisfied. This implies the existence of the matrices $(K,L,M)$ as defined in items \eqref{a}, \eqref{b} and \eqref{c}. We will now show that these matrices indeed constitute a stabilizing controller. Note that by condition \eqref{1}, $\Sigma_{\iso}=\{(A_s,B_s,C_s,D_s)\}$ with
	\begin{equation}
	\label{dataABCD}
	\begin{bmatrix}
	A_s & B_s \\ C_s & D_s
	\end{bmatrix} = \begin{bmatrix}
	X_+ V_1 & X_+ V_2 \\ Y_- V_1 & Y_- V_2
	\end{bmatrix}.
	\end{equation}
	By definition of $K$, $L$ and $M$, the matrices $A_s+B_sM$ and $A_s-LC_s$ are stable and $K = A_s+B_sM-LC_s-LD_sM$. This implies that \eqref{closedloopKLM} is stable since the matrices
	$$
	\begin{bmatrix}
	A_s & B_sM \\ LC_s & A_s+B_sM-LC_s
	\end{bmatrix} \text{ and }\begin{bmatrix}
	A_s + B_sM & B_sM \\ 0 & A_s-LC_s
	\end{bmatrix}
	$$
	are similar \cite[Section 3.12]{Trentelman2001}. We conclude that $(U_-,X,Y_-)$ are informative for stabilization by dynamic measurement feedback and that the recipe given by \eqref{a}, \eqref{b} and \eqref{c} leads to a stabilizing controller $(K,L,M)$.
	
It remains to prove the `only if' part. To this end, suppose that the data $(U_-,X,Y_-)$ are informative for stabilization by dynamic measurement feedback. Let $(K,L,M)$ be such that $\Sigma_{\iso} \subseteq \Sigma_{K,L,M}$. This means that
	$$
	\begin{bmatrix}
	A & BM \\ LC & K+LDM
	\end{bmatrix}
	$$
	is stable for all $(A,B,C,D) \in \Sigma_{\iso}$. Let $\zeta\in\R^n$ and $\eta\in\R^m$ be such that
$$
\bbm\zeta^\top&\eta^\top\ebm\bbm X_- \\ U_-\ebm=0.
$$	
Note that $(A+\zeta\zeta^\top,B+\zeta\eta^\top,C,D)\in \Sigma_{\iso}$ if $(A,B,C,D) \in \Sigma_{\iso}$. Therefore, the matrix $$
	\begin{bmatrix}
	A & BM \\ LC & K+LDM
	\end{bmatrix} + \alpha \begin{bmatrix}
	\zeta\zeta^\top & \zeta\eta^\top M \\ 0 & 0
	\end{bmatrix}
$$
is stable for all $\alpha \in \mathbb{R}$. We conclude that the spectral radius of the matrix 
	\begin{equation*}
	W_\alpha := \frac{1}{\alpha} \begin{bmatrix}
	A & BM \\ LC & K+LDM
	\end{bmatrix} +  \begin{bmatrix}
	\zeta\zeta^\top & \zeta\eta^\top M \\ 0 & 0
	\end{bmatrix}
\end{equation*}
is smaller than $1/\alpha$. By taking the limit as $\alpha \to \infty$, we see that the spectral radius of $\zeta\zeta^\top$ must be zero due to the continuity of spectral radius. Therefore, $\zeta$ must be zero. Since $U_-$ has full column rank, we can conclude that $\eta$ must be zero too. This proves that
condition \eqref{1} and therefore $\Sigma_{\iso} = \{(A_s,B_s,C_s,D_s)\}$. Since the controller $(K,L,M)$ stabilizes $(A_s,B_s,C_s,D_s)$, the pair $(A_s,B_s)$ is stabilizable and $(C_s,A_s)$ is detectable. By \eqref{dataABCD} we conclude that condition \eqref{2} is also satisfied. This proves the theorem.\EP

The following corollary follows from Lemma~\ref{l:equiv} and Theorem~\ref{t:stabUXY} and gives necessary and sufficient conditions for informativity for stabilization by dynamic measurement feedback. Note that we do not make any a priori assumptions on the rank of $U_-$. 

\begin{cor}
	\label{c:stabUXY}
	Let $S$ be any full column rank matrix such that $U_- = S \hat{U}_-$ with $\hat{U}_-$ full row rank $k$. The data $(U_-,X,Y_-)$ are informative for stabilization by dynamic measurement feedback if and only if the following two conditions are satisfied:
	\begin{enumerate}[(i)]
		\item We have
		\[ \rank\begin{bmatrix} X_- \\ \hat{U}_- \end{bmatrix} = n+k.\] 
		Equivalently, there exists a matrix $\begin{bmatrix} V_1 & V_2	\end{bmatrix}$ such that
		\[\begin{bmatrix}
		X_- \\ \hat{U}_-
		\end{bmatrix} \begin{bmatrix}
		V_1 & V_2
		\end{bmatrix} = \begin{bmatrix}
		I & 0 \\ 0 & I
		\end{bmatrix}.\]
		\item The pair $(X_+ V_1, X_+ V_2)$ is stabilizable and $(Y_- V_1, X_+ V_1)$ is detectable. 
	\end{enumerate}
	Moreover, if the above conditions are satisfied, a stabilizing controller $(K,L,M)$ is constructed as follows:
	\begin{enumerate}[(a)]
		\item Select a matrix $\hat{M}$ such that $X_+ (V_1 + V_2 \hat{M})$ is stable. Define $M := S \hat{M}$. \label{c:a}
		\item Choose a matrix $L$ such that $(X_+ - L Y_-) V_1$ is stable. \label{c:b}
		\item Define $K := (X_+ - LY_-)(V_1 + V_2 \hat{M})$. \label{c:c}
	\end{enumerate} 
\end{cor}

\begin{remark}
	In the previous corollary it is clear that the system matrices of the data-generating system are related to the data via
	\[ \begin{bmatrix} A_s & B_s S \\ C_s & D_s S \end{bmatrix} = \begin{bmatrix} X_+ \\ Y_- \end{bmatrix} \begin{bmatrix} V_1 & V_2\end{bmatrix}.\] 
	Therefore the corollary shows that informativity for stabilization by dynamic measurement feedback requires that $A_s$ and $C_s$ can be identified uniquely from the data. However, this does not hold for $B_s$ and $D_s$ in general.
\end{remark}

\subsection{Stabilization using input and output data}

Recall that we consider a system of the form \eqref{inoutsys}. When given input, state and output data, any system $(A,B,C,D)$ consistent with these data satisfies 
\begin{equation}\label{eq:iso}
\begin{bmatrix} X_+ \\ Y_- \end{bmatrix} = \begin{bmatrix} A & B \\ C & D\end{bmatrix} \begin{bmatrix} X_- \\ U_- \end{bmatrix}.
\end{equation}
In this section, we will consider the situation where we have access to input and output measurements only. Moreover, we assume that the data are collected on a single time interval, i.e. $q=1$. This means that our data are of the form $(U_-,Y_-)$, where
\begin{subequations}\label{eq:uy data}
	\begin{align}
	U_- &:= \bbm u(0) & u(1) & \cdots & u(T-1) \ebm \\
	Y_- &:= \bbm y(0) & y(1) & \cdots & y(T-1) \ebm . 
	\end{align}
\end{subequations}

Again, we are interested in informativity of the data, this time given by $(U_-,Y_-)$. Therefore we wish to consider the set of all systems of the form \eqref{inoutsys} with the state space dimension\footnote{The state space dimension of the system may be known a priori. In the case that it is not, it can be computed using subspace identification methods, see e.g. \cite[Theorem 2]{vanOverschee1996}.} $n$ that admit the same input/output data. This leads to the following set of consistent systems:
\[ \Sigma_{\io} := \left\{ (A,B,C,D) \mid \exists X\in\mathbb{R}^{n\times (T+1)} \textrm{ s.t. } \eqref{eq:iso} \textrm{ holds} \right\}. \]
As in the previous section, we wish to find a controller of the form \eqref{compensator} that stabilizes the system. This means that, in line with Definition~\ref{def:par informativity}, we have the following notion of informativity:

\begin{definition}
	We say the data $(U_-,Y_-)$ are \emph{informative for stabilization by dynamic measurement feedback} if there exist matrices $K$, $L$ and $M$ such that $\Sigma_{\io} \subseteq \Sigma_{K,L,M}$.
\end{definition}

%Unlike in the input/state case, $\Sigma_{\io}$ is never a singleton. Clearly, if $(A,B,C,D)\in\Sigma_{\io}$ with state sequence $X$, then for any invertible matrix $S$, the system $(SAS^{-1}, SB,CS^{-1},D)$ is an element of $\Sigma_{\io}$ with state sequence $SX$. However, if $(K,L,M)$ is a stabilizing controller for $(A,B,C,D)$, it also stabilizes all systems similar to $(A,B,C,D)$.  

In order to obtain conditions under which $(U_-,Y_-)$ are informative for stabilization, it may be tempting to follow the same steps as in Section~\ref{sec:iso measurements}. In that section we first proved that we can assume without loss of generality that $U_-$ has full row rank. Subsequently, Theorem~\ref{t:stabUXY} and Corollary~\ref{c:stabUXY} characterize informativity for stabilization by dynamic measurement feedback based on input, state and output data. It turns out that we can perform the first of these two steps for input/output data as well. Indeed, in line with Lemma~\ref{l:equiv}, we can state the following: 

\begin{lemma}
	\label{l:equiv2}
	Consider the data $(U_-,Y_-)$ and the corresponding set $\Sigma_{\io}$. Let $S$ be a matrix of full column rank such that $U_- = S \hat{U}_-$ with $\hat{U}_-$ a matrix of full row rank. 
	
	Then the data $(U_-,Y_-)$ are informative for stabilization by dynamic measurement feedback if and only if the data $(\hat{U}_-,Y_-)$ are informative for stabilization by dynamic measurement feedback. 
\end{lemma}
The proof of this lemma is analogous to that of Lemma~\ref{l:equiv} and therefore omitted. Lemma~\ref{l:equiv2} implies that without loss of generality we can consider the case where $U_-$ has full row rank. 

In contrast to the first step, the second step in Section~\ref{sec:iso measurements} relies heavily on the affine structure of the considered set $\Sigma_{\iso}$. Indeed, the proof of Theorem~\ref{t:stabUXY} makes use of the fact that $\Sigma^0_{\iso}$ is a subspace. However, the set $\Sigma_{\io}$ is not an affine set. This means that it is not straightforward to extend the results of Corollary~\ref{c:stabUXY} to the case of input/output measurements. 

Nonetheless, under certain conditions on the input/output data it is possible to construct the corresponding state sequence $X$ of \eqref{inoutsys} up to similarity transformation. In fact, state reconstruction is one of the main themes of subspace identification, see e.g. \cite{Moonen1989,vanOverschee1996}. The construction of a state sequence would allow us to reduce the problem of stabilization using input/output data to that with input, state and output data. The following result gives sufficient conditions on the data $(U_-,Y_-)$ for state construction. 

To state the result, we will first require a few standard pieces of notation. First, let $f(0),\ldots, f(T-1)$ be a signal and $\ell <T$, then we define the \textit{Hankel matrix of depth} $\ell$ as
\[ \calH_\ell(f)= \begin{bmatrix} f(0) & f(1) & \cdots & f(T-\ell) \\ f(1) & f(2) & \cdots & f(T-\ell+1) \\ \vdots & \vdots &  &\vdots \\ f(\ell-1) & f(\ell) & \cdots &f(T-1) \end{bmatrix}. \]
Given input and output data of the form \eqref{eq:uy data}, and $k$ such that $2k<T$ we consider $\calH_{2k}(u)$ and $\calH_{2k}(y)$. Next, we partition our data into so-called `\textit{past}' and `\textit{future}' data as
\[ \calH_{2k}(u) = \begin{bmatrix} U_p \\U_f \end{bmatrix}, \quad  \calH_{2k}(y) = \begin{bmatrix} Y_p \\ Y_f \end{bmatrix},\]
where $U_p,U_f,Y_p$ and $Y_f$ all have $k$ block rows. Let $x(0),\ldots, x(T)$ denote the state trajectory of \eqref{inoutsys} compatible with a given $(U_-,Y_-)$. We now denote
\begin{align*}
X_p &= \begin{bmatrix} x(0) & \cdots & x(T-2k) \end{bmatrix}, \\  
X_f &= \begin{bmatrix} x(k) & \cdots & x(T-k)\end{bmatrix}.
\end{align*}

Lastly, let $\rs(M)$ denote the row space of the matrix $M$. Now we have the following result, which is a rephrasing of \cite[Theorem 3]{Moonen1989}.

\begin{theorem}
	Consider the system \eqref{inoutsys} and assume it is minimal. Let the input/output data $(U_-,Y_-)$ be as in \eqref{eq:uy data}. Assume that $k$ is such that $n<k<\tfrac{1}{2}T$. If
	\begin{equation}\label{eq:rank cond} \rank \begin{bmatrix} \calH_{2k}(u) \\\calH_{2k}(y)\end{bmatrix} = 2km+n,\end{equation} 
	then 
	\[ \rs (X_f) = \rs \left(\begin{bmatrix} U_p \\ Y_p\end{bmatrix}\right) \cap \rs \left(\begin{bmatrix} U_f \\ Y_f\end{bmatrix}\right), \]
	and this row space is of dimension $n$. 
\end{theorem}
Under the conditions of this theorem, we can now find the `true' state sequence $X_f$ up to similarity transformation. That is, we can find $\bar{X} = SX_f$ for some unknown invertible matrix $S$. This means that, under these conditions, we obtain an input/state/output trajectory given by the matrices
\begin{subequations}\label{eq:uyx bar data}
\begin{align}
\bar{U}_- &= \bbm u(k) & u(k+1) & \cdots & u(T-k-1) \ebm, \\
\bar{Y}_- &= \bbm y(k) & y(k+1) & \cdots & y(T-k-1) \ebm, \\
\bar{X} &= S\bbm x(k) & x(k+1) & \cdots & x(T-k)\ebm.
\end{align}
\end{subequations}

We can now state the following sufficient condition for informativity for stabilization with input/output data. 

\begin{cor}\label{c:inf uy}
	Consider the system \eqref{inoutsys} and assume it is minimal. Let the input/output data $(U_-,Y_-)$ be as in \eqref{eq:uy data}. Assume that $k$ is such that $n<k<\tfrac{1}{2}T$. Then the data $(U_-,Y_-)$ are informative for stabilization by dynamic measurement feedback if the following two conditions are satisfied:
	\begin{enumerate}[(i)] 
		\item The rank condition \eqref{eq:rank cond} holds.
		\item The data $(\bar{U}_-,\bar{X},\bar{Y}_-)$, as defined in \eqref{eq:uyx bar data}, are informative for stabilization by dynamic measurement feedback. 
	\end{enumerate}
	Moreover, if these conditions are satisfied, a stabilizing controller $(K,L,M)$ such that $\Sigma_{\io}\subseteq \Sigma_{K,L,M}$ can be found by applying Corollary~\ref{c:stabUXY} \eqref{c:a},\eqref{c:b},\eqref{c:c} to the data $(\bar{U}_-,\bar{X},\bar{Y}_-)$.	
\end{cor}
The conditions provided in Corollary~\ref{c:inf uy} are sufficient, but not necessary for informativity for stabilization by dynamic measurement feedback. In addition, it can be shown that data satisfying these conditions are also informative for system identification, in the sense that $\Sigma_{\io}$ contains only the `true' system \eqref{inoutsys} and all systems similar to it. 

An interesting question is whether the conditions of Corollary~\ref{c:inf uy} can be sharpened to necessary and sufficient conditions. In this case it would be of interest to investigate whether such conditions are weaker than those for informativity for system identification. 

At this moment, we do not have a conclusive answer to the above question. However, we note that even for subspace identification there are no known necessary and sufficient conditions for data to be informative, although several sufficient conditions exist, e.g. \cite[Theorems 3 and 5]{Moonen1989}, \cite[Theorem 2]{vanOverschee1996} and  \cite[Theorems 3 and 4]{Verhaegen1992}.

\section{Conclusions and future work}\label{sec:conc}
Results in data-driven control should clearly highlight the differences and possible advantages as compared to system identification paired with model-based control. One clear advantage of data-driven control is its capability of solving problems in the presence of data that are not informative for system identification. Therefore, informativity is a very important concept for data-driven analysis and control. 

In this paper we have introduced a comprehensive framework for studying informativity problems. We have applied this framework to analyze several system-theoretic properties on the basis of data. The same framework was used to solve multiple data-driven control problems. 

After solving these problems, we have made the comparison between our data-driven methods, and the `classical' combination of identification and model-based control. We have shown that for many analysis and control problems, such as controllability analysis and stabilization, the data-driven approach can indeed be performed on data that are not informative for system identification. On the other hand, for data-driven linear quadratic regulation it has been shown that informativity for system identification is a necessary condition. This effectively means that for this data-driven control problem, we have given a theoretic justification for the use of persistently exciting data.

\subsection*{Future work} 
Due to the generality of the introduced framework, many different problems can be studied in a similar fashion: one could consider different types of data, where more results based on only input and output data would be particularly interesting. Many other system-theoretic properties could be considered as well, for example, analyzing passivity or tackling robust control problems based on data. 

It would also be of interest to generalize the model class under consideration. One could, for instance, consider larger classes of systems like differential algebraic or polynomial systems. On the other hand, the class under consideration can also be made smaller by prior knowledge of the system. For example, the system might have an observed network structure, or could in general be parametrized. 

A framework similar to ours could be employed in the presence of disturbances, which is a problem of practical interest. A study of data-driven control problems in this situation is particularly interesting, because system identification is less straightforward. 
We note that data-driven stabilization under measurement noise has been studied in \cite{DePersis2019} and under unknown disturbances in \cite{Berberich2019c}. Additionally, the data-driven LQR problem is popular in the machine learning community, where it is typically assumed that the system is influenced by (Gaussian) process noise, see e.g. \cite{Dean2019}.

In this paper, we have assumed that the data are given. Yet another problem of practical interest is that of \textit{experiment design}, where inputs need to be chosen such that the resulting data are informative. In system identification, this problem led to the notion of persistence of excitation. For example, it is shown in \cite{Willems2005} that the rank condition \eqref{eq:inf for sys ident} can be imposed by injecting an input sequence that is persistently exciting of order $n+1$. However, as we have shown, this rank condition is not necessary for some data-driven control problems, like stabilization by state feedback. The question therefore arises whether we can find tailor-made conditions on the input only, that guarantee informativity for data-driven control.

\section*{Acknowledgements}
We thank our colleague Prof. Arjan van der Schaft for his valuable comments.
%The informativity problem is a generalization of a problem from the field of system identification. In this field, uniquely identifying a system from data is often built on two central concepts: \textit{Identifiability of the system} and \textit{informativity of the data}. The first of these is the problem of finding conditions on the system under which there exists data for which the system can be uniquely identified. In terms of the framework of this paper, the \textit{identifiability problem} can be stated as: 

%\begin{problem}[Identifiability problem]
%	Provide necessary and sufficient conditions on the \textit{system} under which there exists data $\calD$ that are informative for property $\calP$. 
%\end{problem}

%Especially in cases where conditions for identifiability for system identification are not known, the study of identifiability might need to precede the question of informativity. 

\bibliographystyle{IEEEtran}
\bibliography{references}

\end{document}